\documentclass[11pt]{article}

\usepackage{amsfonts}
\usepackage{amscd}
\usepackage{amssymb}
\usepackage{amsthm}
\usepackage{amsmath}

\input prepictex
\input pictex
\input postpictex

\theoremstyle{plain}
\newtheorem{thm}{Theorem}[section]

\newtheorem{prop}[thm]{Proposition}

\theoremstyle{definition}

\newtheorem{remark}[thm]{Remark}
\theoremstyle{example}

\theoremstyle{remark}

\numberwithin{equation}{section}

\def\cB{\mathcal{B}}

\def\cP{\mathcal{P}}

\def\CC{\mathbb{C}}

\def\RR{\mathbb{R}}

\def\ZZ{\mathbb{Z}}

\def\fh{\mathfrak{h}}

\def\Ind{\mathrm{Ind}}

\def\wt{\mathrm{wt}}

\makeatletter
\renewcommand{\@makefnmark}{\mbox{\textsuperscript{}}}
\makeatother

\title{A combinatorial formula for Macdonald polynomials}
\author{Arun Ram \\
Department of Mathematics\\
University of Wisconsin\\
Madison, WI 53706 \\
ram@math.wisc.edu \\
and \\
Department of Mathematics and Statistics \\
University of Melbourne \\
Parkville VIC 3010 Australia
\and
Martha Yip \\
Department of Mathematics\\
University of Wisconsin\\
Madison, WI 53706 \\
yip@math.wisc.edu}
\date{}

\begin{document}

\maketitle

$$\textit{Dedicated to Adriano Garsia}$$

\bigskip

\begin{abstract}
Abstract.  In this paper we use the combinatorics of alcove walks to give a uniform
combinatorial formula for Macdonald polynomials for all Lie types.  These formulas are
generalizations of the formulas of Haglund-Haiman-Loehr for Macdonald polynoimals
of type $GL_n$.  At $q=0$ these formulas specialize to the formula of Schwer
for the Macdonald spherical function in terms of positively folded alcove walks and
at $q=t=0$ these
formulas specialize to the formula for the Weyl character in terms of the Littelmann path model
(in the positively folded gallery form of Gaussent-Littelmann).
\end{abstract}


\begin{section}{Introduction}

The Macdonald polynomials were introduced in the mid 1980s \cite{Mac1} \cite{Mac2}
as a remarkable
family of orthogonal polynomials generalizing the spherical functions for a $p$-adic
group, the Weyl characters, the Jack polynomials and the zonal polynomials.  In
the early 1990s Cherednik \cite{Ch} introduced the double affine Hecke algebra (the DAHA) and
used it as a tool to prove conjectures of Macdonald.
The DAHA is a fundamental tool for studying Macdonald polynomials.
Using the DAHA,
the nonsymmetric Macdonald polynomials $E_\mu$ can be constructed by
applying products of ``intertwining operators'' $\tau_i^\vee$ to the generator $\mathbf{1}$ of the
polynomial representation of the DAHA  (see \cite[Prop.\ 6.13]{Hai}),
and the symmetric Macdonald polynomials $P_\mu$ can
then be constructed from the $E_\mu$ by ``symmetrizing'' (see \cite[Remarks after (6.8)]{Mac3}).

Of recent note in the theory of Macdonald polynomials has been the success of
Haglund-Haiman-Loehr in giving, in the type $GL_n$ case,
explicit combinatorial formulas for the expansion of Macdonald
polynomials in terms of monomials.  These formulas were conjectured by
J.\ Haglund and proved by Haglund-Haiman-Loehr in \cite{HHL1} and \cite{HHL2}.
The papers \cite{GR} and \cite{Hai} are excellent survey articles discussing these
developments.

Following a key idea of C.\ Schwer \cite{Sc}, the paper \cite{Ra} developed a
combinatorics for working in the affine Hecke algebra, the \emph{alcove walk model}.
It turns out that this combinatorics is the ideal tool for expansion of products of
intertwining operators
in the DAHA. These expansions, when applied to the generator
of the polynomial representation of the DAHA, give formulas for the Macdonald polynomials
which are generalizations, to all root systems,
 of the formulas obtained by Haglund-Haiman-Loehr \cite{HHL1} \cite{HHL2}
in type $GL_n$.

At $q=0$ the symmetric Macdonald polynomials are the \emph{Hall-Littlewood
polynomials} or the \emph{Macdonald spherical functions}. These
are the spherical functions for $G/K$, where $G$ is a $p-$adic group and $K$ is a
maximal compact subgroup.  The work of Schwer \cite[Thm.\ 1.1]{Sc}
provided fomulas for the expansion of the Macdonald spherical functions in terms
of positively folded alcove walks.  See \cite[Thm.\ 4.2(a)]{Ra} for
a description of the Schwer-KLM formula in terms of the alcove walk model.
The formula
for Macdonald polynomials which we give in Theorem \ref{Pexpansion}
reduces to the Schwer formula at $q=0$.

At $q=t=0$ the symmetric Macdonald polynomials
are the \emph{Weyl characters} or \emph{Schur functions}.
In this case our formula for the Macdonald polynomial
specializes to the formula for the Weyl character in terms of the Littelmann path model
(in the maximal dimensional positively folded gallery form of Gaussent-Littelmann
\cite[Cor.\ 1 p.\ 62]{GL}).

It is interesting to note that, in the formulas for the symmetric Macdonald polynomials,
the negative folds and the positive folds play an equal role.  It is known \cite{GL} that the
alcove walks with only positive folds contain detailed information about the geometry
of Mirkovi\'c-Vilonen intersections in the loop Grassmannian.  It is tantalizing to wonder
whether the alcove walks with both positive and negative folds play a similar role in the geometry of
flag varieties for reductive groups over two dimensional local fields and whether the expansions
of Macdonald polynomials in this paper are shadows of geometric decompositions.

The papers \cite{GL} and \cite {Ra} explain how the combinatorics of alcove walks is
almost equivalent to the combinatorics of crystal bases and Kashiwara operators (at least
for the positively folded alcove walks of maximal dimension).  Our expansions of Macdonald
polynomials in terms of alcove walks give insight into possible relationships between
Macdonald polynomials and crystal and canonical bases.

This research was partially supported by the National
Science Foundation (NSF) under grant DMS-0353038 at the University of Wisconsin, Madison.
We thank
the NSF for continuing support of our research.  This paper was completed while the authors
were in residence at  the
special semester in Combinatorial Representation Theory at Mathematical Sciences
Research Institute (MSRI).  It is a pleasure to thank MSRI for hospitality, support and a wonderful
and stimulating working environment.
A.\ Ram thanks S.\ Griffeth for many many instructive conversations about
double affine Hecke algebras and Macdonald polynomials, without which this paper would
never have been possible.

\end{section}

\begin{section}{Double Weyl groups, braid groups and Hecke algebras}

In this section we review the basic definitions and notations for affine Weyl groups and
double affine Hecke algebras
following the expositions in \cite{Ra}, \cite{Ch}, \cite{Mac4} and \cite {Hai}.  Following the
definitions we prove Theorem \ref{mainthm},
a formula for the expansion of products of intertwining operators in the DAHA.  This
formula is a ``lift into the DAHA'' of the expansions of Macdonald polynomials given in
Section 3.

\subsection{Double affine Weyl groups}

Let $\fh_\ZZ$ be a $\ZZ$-lattice with an action of a finite subgroup
$W_0$ of $GL(\fh_\ZZ)$ generated by reflections. Then $W_0$ acts on $\fh_\ZZ^*$ by
\begin{equation}
\langle w\mu, \lambda^\vee\rangle = \langle \mu, w^{-1}\lambda^\vee\rangle,
\qquad\hbox{where}\quad
\langle \lambda^\vee, \mu\rangle = \mu(\lambda^\vee)
\ \ \hbox{for $\lambda^\vee\in \fh_\ZZ$, $\mu\in \fh_\ZZ^*$.}
\end{equation}
Let $R^+\subseteq \fh_\ZZ^*$ and  $(R^\vee)^+\subseteq \fh_\ZZ$ denote fixed choices
of the positive roots and
the positive coroots so that the reflections
$s_\alpha$ in $W_0$ act on $\fh_\ZZ$ and
on $\fh_\ZZ^*$ by
\begin{equation}
s_\alpha \lambda = \lambda - \langle\lambda,\alpha^\vee\rangle\alpha
\qquad\hbox{and}\qquad
s_\alpha \lambda^\vee = \lambda^\vee - \langle \lambda^\vee,\alpha\rangle \alpha^\vee,
\qquad\hbox{respectively.}
\end{equation}

The groups
\begin{equation}
X = \{ X^\mu\ |\ \mu\in \fh_\ZZ^*\}
\quad\hbox{and}\quad Y = \{ Y^{\lambda^\vee}\ |\ \lambda^\vee\in
\fh_\ZZ\}
\end{equation}
with
\begin{equation}\label{XYrels}
X^{\mu}X^{\nu} = X^{\mu+\nu}
\quad\hbox{and}\quad
Y^{\lambda^\vee}Y^{\sigma^\vee} = Y^{\lambda^\vee+\sigma^\vee}
\end{equation}
are the groups $\fh^*_\ZZ$ and $\fh_\ZZ$ respectively,
except written multiplicatively, and the semidirect product
\begin{equation}
W_0\ltimes (X\times Y) = \{ X^\mu w Y^{\lambda^\vee}
\ |\ w\in W_0, \mu\in \fh_\ZZ^*, \lambda^\vee\in \fh_\ZZ\}
\end{equation}
has additional relations
\begin{equation}\label{WXYcomm}
\qquad
wX^\mu = X^{w\mu}w
\quad\hbox{and}\quad
wY^{\lambda^\vee} = Y^{w\lambda^\vee}w,
\end{equation}
for $w\in W_0$, $\mu\in \fh_\ZZ^*$ and $\lambda^\vee\in \fh_\ZZ$.

Assume that the action of $W_0$ on $\fh_\CC = \CC\otimes_\ZZ \fh_\ZZ$ is
irreducible.  The \emph{double affine Weyl group} $\widetilde{W}$ is the universal central extension
of $W_0\ltimes(X\times Y)$. If $e$ is the smallest integer such that
$\langle \lambda^\vee,\mu\rangle\in \frac{1}{e}\ZZ$ for all $\lambda^\vee\in \fh_\ZZ$
and $\mu\in \fh_\ZZ^*$ then $\widetilde{W}$ is presented by
$$\widetilde{W} = \{ q^k  X^\mu w Y^{\lambda^\vee}\ |\
k\in \hbox{$\frac{1}{e}$}\ZZ, \mu\in \fh_\ZZ^*, \lambda^\vee\in \fh_\ZZ, w\in W_0\}$$
with \eqref{XYrels}, \eqref{WXYcomm} and
\begin{equation}\label{dblW}
X^\mu Y^{\lambda^\vee}
= q^{\langle\lambda^\vee,\mu\rangle}Y^{\lambda^\vee}X^\mu,
\qquad\hbox{for $\mu\in \fh_\ZZ^*$, $\lambda^\vee\in \fh_\ZZ$.}
\end{equation}
The subgroup
$\{q^kX^\mu Y^{\lambda^\vee}\ |\ k\in \frac{1}{e}\ZZ,\, \mu\in \fh_\ZZ,\,
\lambda^\vee\in \fh_\ZZ\}$ is a Heisenberg group and
\begin{equation}
W = \{ X^\mu w\ |\ \mu\in \fh_\ZZ, w\in W_0\}
\qquad\hbox{and}\quad
W^\vee = \{ w Y^{\lambda^\vee}\ |\ \lambda^\vee\in \fh_\ZZ, w\in W_0\}
\end{equation}
are \emph{affine Weyl groups} inside $\widetilde{W}$.  Letting
\begin{equation}\label{qdeltad}
q = X^\delta = Y^{-d}
\end{equation}
and extending the notation of \eqref{WXYcomm} gives actions of
$W^\vee$ on $\fh_\ZZ^*+\ZZ \delta$ and $W$ on $\fh_\ZZ\oplus \ZZ d$ with
\begin{equation}\label{Waction}
Y^{\lambda^\vee} \mu = \mu - \langle \mu, \lambda^\vee\rangle \delta
\qquad\hbox{and}\qquad
X^\mu \lambda^\vee = \lambda^\vee - \langle \lambda^\vee,\mu\rangle d.
\end{equation}
Let $\varphi\in R$ be the highest root and $\varphi^\vee\in R^\vee$ the highest coroot and let
\begin{equation}\label{s0defn}
s_0 = Y^{\varphi^\vee} s_\varphi
\qquad\hbox{and}\qquad
s_0^\vee = X^{\varphi}s_{\varphi^\vee}.
\end{equation}
Let
\begin{equation}\label{a0defn}
\alpha_0 = -\varphi+\delta,
\quad
\alpha_0^\vee = -\varphi^\vee + d,
\qquad
\langle d,\mu\rangle=0,
\quad
\langle \lambda^\vee, \delta\rangle=0,
\qquad
\langle d, \delta\rangle = 0,
\end{equation}
so that
\begin{equation}\label{s0action}
s_0\mu = \mu - \langle \mu, \alpha_0^\vee\rangle\alpha_0
\qquad\hbox{and}\qquad
s_0^\vee\lambda^\vee = \lambda^\vee - \langle \lambda^\vee, \alpha_0\rangle\alpha^\vee_0.
\end{equation}

The \emph{alcoves} of $\fh_\RR^*=\RR\otimes_\ZZ \fh^*_\ZZ$
are the connected components of
\begin{equation}
\fh_\RR^*\backslash \left( \bigcup_{\alpha^\vee\in (R^\vee)^+,\, j\in \ZZ}
\fh^{\alpha^\vee+jd}\right)
\qquad\hbox{where}\qquad
\fh^{\alpha^\vee+jd}=\{  x\in \fh_\RR^*\ |\ \langle x, \alpha^\vee\rangle = -j \}.
\end{equation}
The action of $W=\{X^\mu w\ |\ \mu\in \fh_\ZZ^*, w\in W_0\}$ on $\fh_\RR^*$ given by
\begin{equation}\label{translaction}
X^\mu \cdot \nu = \nu+\mu
\qquad\hbox{and}\qquad
w\cdot \nu = w\nu,
\qquad\hbox{for $w\in W_0$, $\mu\in \fh^*_\ZZ$ and $\nu\in \fh^*_\RR$,}
\end{equation}
sends alcoves to alcoves;
$s^\vee_0,\ldots, s^\vee_n$ are the reflections in the walls $\fh^{\alpha_0^\vee}, \ldots,
\fh^{\alpha_n^\vee}$ of the fundamental alcove
\begin{align}
1 &= \{ x\in \fh_\RR^*\ |\ \hbox{$\langle x,\alpha^\vee_i\rangle \ge 0$, for $i=0,1,\ldots, n$}\};
\qquad\hbox{and} \\
\ell(v) &= \hbox{(number of hyperplanes between $1$ and $v$)}
\end{align}
is the \emph{length} of $v\in W$.
Let $\Omega^\vee$ be the set of length zero elements of $W$.
The affine Weyl group $W$ has an alternate presentation by generators
$s^\vee_0, s^\vee_1,\ldots, s^\vee_n$ and $\Omega^\vee$ with relations
\begin{equation}
(s^\vee_i)^2 = 1,
\qquad \underbrace{s^\vee_is^\vee_j\cdots}_{m^\vee_{ij}} =
\underbrace{s^\vee_js^\vee_i\cdots}_{m^\vee_{ij}},
\qquad\hbox{and}\qquad
g^\vee s^\vee_i(g^\vee)^{-1}=s^\vee_{\sigma^\vee(i)},\quad\hbox{for $g^\vee\in\Omega^\vee$},
\end{equation}
where $\pi/m^\vee_{ij}$ is the angle between $\fh^{\alpha_i^\vee}$ and $\fh^{\alpha_j^\vee}$ and
$\sigma^\vee$ denotes the permutation of the $\fh^{\alpha_i^\vee}$ induced by the action of
$g^\vee$.
If $\Omega^\vee\times \fh_\RR^*$ is $|\Omega^\vee|$ copies of $\fh_\RR^*$ (sheets), with
$\Omega^\vee$ acting by switching sheets then there is a bijection
\begin{equation}\label{Wtoalcoves}
W \longleftrightarrow \{ \hbox{alcoves in $\Omega^\vee\times \fh_\RR^*$}\}
\end{equation}
and we will often identify $v\in W$ with the corresponding alcove in $\Omega^\vee\times \fh_\RR^*$.
The pictures illustrating this bijection in type $SL_3$ are displayed in the appendix.

The \emph{periodic orientation} is the orientation of the hyperplanes $\fh^{\alpha^\vee+kd}$
such that
\begin{equation}\label{perorient}
\begin{array}{l}
\hbox{(a) $1$ is on the positive side of $\fh^{\alpha^\vee}$ for $\alpha^\vee\in (R^\vee)^+$,}
\qquad\qquad\qquad\qquad\qquad \\ \\
\hbox{(b) $\fh^{\alpha^\vee+kd}$ and $\fh^{\alpha^\vee}$ have parallel orientations.}
\end{array}
\end{equation}
The pictures in the appendix illustrate the periodic orientation for type $SL_3$.

A similar ``pictorial''
viewpoint applies to the group $W^\vee$ acting on $\Omega\times \fh_\RR$ where
$\fh_\RR = \RR\otimes_\ZZ \fh_\ZZ$
and $\Omega$ is the set of length zero elements of $W^\vee$.
Then $W^\vee$ has an alternate presentation by generators
$s_0, s_1,\ldots, s_n$ and $\Omega$ with relations
\begin{equation}
s_i^2 = 1,
\qquad
\underbrace{s_is_j\cdots}_{m_{ij}} =
\underbrace{s_js_i\cdots}_{m_{ij}},
\qquad\hbox{and}\qquad
g s_ig^{-1}=s_{\sigma(i)},\quad\hbox{for $g\in\Omega$},
\end{equation}
where $\pi/m_{ij}$ is the angle between $\fh^{\alpha_i}$ and $\fh^{\alpha_j}$
and $\sigma$ denotes the permutation of the $\fh^{\alpha_i}$ induced by the action
of $g$.

\subsection{Double affine braid groups}

The \emph{double affine braid group} $\tilde \cB$ is the group generated by
$T_0,\ldots, T_n$, $\Omega$ and $X$ with relations
\begin{equation}\label{dblbraidreln2}
\underbrace{T_iT_j\cdots}_{m_{ij}} =
\underbrace{T_jT_i\cdots}_{m_{ij}}, \quad g T_ig^{-1}=T_{\sigma(i)},
\quad g X^\mu = X^{g\mu}g,
\end{equation}
for $g\in\Omega$, and
\begin{equation}\label{dblbraidreln1}
\begin{array}{cl}
T_iX^\mu = X^{s_i\mu}T_i, &\textrm{if } \langle
\mu,\alpha_i^\vee\rangle = 0, \\
T_iX^\mu T_i = X^{s_i\mu}, &\textrm{if } \langle
\mu,\alpha_i^\vee\rangle = 1,
\end{array}
\qquad\hbox{for $i=0,1,\ldots, n$,}
\end{equation}
where the action of $W^\vee$ on $\fh^*_\ZZ\oplus \ZZ\delta$ is as
in \eqref{Waction}.
The element
\begin{equation}
q = X^\delta
\quad\hbox{is in the center of $\tilde \cB$.}
\end{equation}

For $w\in W^\vee$, view a reduced word $w =
gs_{i_1}\cdots s_{i_\ell}$ as a minimal length path $p$ from the
fundamental alcove to $w$ in $\fh_\RR$ and define
\begin{equation}
Y^w = g(T_{i_1})^{\epsilon_1}\cdots (T_{i_\ell})^{\epsilon_\ell},
\qquad\hbox{with}\qquad
\epsilon_k = \begin{cases}
+1, &\hbox{if the $k$th step of $p$ is}\ \
\beginpicture
\setcoordinatesystem units <1cm,1cm>                    
\setplotarea x from -0.8 to 0.8, y from -0.5 to 0.5     
\put{$\scriptstyle{-}$}[b] at -0.4 0.25 \put{$\scriptstyle{+}$}[b]
at 0.4 0.25
\plot  0 -0.4  0 0.5 /
\arrow <5pt> [.2,.67] from -0.5 0 to 0.5 0   %
\endpicture
, \\
-1, &\hbox{if the $k$th step of $p$ is}\ \
\beginpicture
\setcoordinatesystem units <1cm,1cm>                    
\setplotarea x from -0.8 to 0.8, y from -0.5 to 0.5     
\put{$\scriptstyle{-}$}[b] at -0.4 0.25 \put{$\scriptstyle{+}$}[b]
at 0.4 0.25
\plot  0 -0.4  0 0.5 /
\arrow <5pt> [.2,.67] from 0.5 0 to -0.5 0   %
\endpicture
,
\end{cases}
\end{equation}
with respect to the \emph{periodic orientation} (see \eqref{perorient} and the pictures in the
appendix).
For $v\in W$, view a reduced word $v =
g^\vee s^\vee_{i_1}\cdots s^\vee_{i_\ell}$ as a minimal length path $p^\vee$ from the
fundamental alcove to $v$ in $\fh_\RR^*$ and define
\begin{equation}
X^{v} = g^\vee(T^\vee_{i_1})^{\epsilon^\vee_1}\cdots (T^\vee_{i_\ell})^{\epsilon^\vee_\ell},
\quad\hbox{with}\quad
\epsilon^\vee_k = \begin{cases}
-1, &\hbox{if the $k$th step of $p^\vee$ is}\ \
\beginpicture
\setcoordinatesystem units <1cm,1cm>                    
\setplotarea x from -0.8 to 0.8, y from -0.5 to 0.5     
\put{$\scriptstyle{-}$}[b] at -0.4 0.25 \put{$\scriptstyle{+}$}[b]
at 0.4 0.25
\plot  0 -0.4  0 0.5 /
\arrow <5pt> [.2,.67] from -0.5 0 to 0.5 0   %
\endpicture
, \\
+1, &\hbox{if the $k$th step of $p^\vee$ is}\ \
\beginpicture
\setcoordinatesystem units <1cm,1cm>                    
\setplotarea x from -0.8 to 0.8, y from -0.5 to 0.5     
\put{$\scriptstyle{-}$}[b] at -0.4 0.25 \put{$\scriptstyle{+}$}[b]
at 0.4 0.25
\plot  0 -0.4  0 0.5 /
\arrow <5pt> [.2,.67] from 0.5 0 to -0.5 0   %
\endpicture
,
\end{cases}
\end{equation}

Let
$T_i^\vee = T_i$, for $i=1,2,\ldots, n$,
\begin{equation}\label{gT0}
g^\vee = X^{\omega_g}T^\vee_{w_gw_0},
\qquad
(T_0^\vee)^{-1} = X^{\varphi}T^\vee_{s_\varphi},
\qquad
g  = Y^{\omega^\vee_g}T^{-1}_{w_0w_g},
\qquad
T_0 = Y^{\varphi^\vee}T^{-1}_{s_\varphi}.
\end{equation}
where $\varphi$ and $\varphi^\vee$ are as in \eqref{s0defn} and,
using the action in \eqref{translaction}, $\omega_g= g^\vee\cdot 0$ and $w_g$ is the longest
element of the stabilizer of $\omega_g$ in $W_0$.

The following theorem, discovered by Cherednik \cite[Thm.\ 2.2]{Ch}, is proved in \cite[3.5-3.7]{Mac4},
in \cite{Io}, and in \cite[4.13-4.18]{Hai}.

\begin{thm}\label{duality} (Duality)  Let $Y^d =q^{-1}$.
The double braid group
 $\tilde \cB$ is generated by $T_0^\vee, T_1^\vee,\ldots, T_n^\vee$, $\Omega^\vee$ and
$Y$ with relations
\begin{equation}\label{dualdblbraidreln}
\underbrace{T^\vee_iT^\vee_j\cdots}_{m^\vee_{ij}} =
\underbrace{T^\vee_jT^\vee_i\cdots}_{m^\vee_{ij}},
\qquad g^\vee T^\vee_i(g^\vee)^{-1}=T^\vee_{\sigma^\vee(i)},
\qquad g^\vee Y^{\lambda^\vee} = Y^{g^\vee\lambda^\vee}g^\vee,
\end{equation}
for $g^\vee\in\Omega^\vee$,
and
\begin{equation}\label{doubleBraidTY}
\begin{array}{cl}
T_i^\vee Y^{\lambda^\vee} = Y^{s^\vee_i\lambda^\vee}T_i^\vee, &\textrm{if } \langle
\lambda^\vee,\alpha_i\rangle = 0, \\
(T_i^\vee)^{-1}Y^{\lambda^\vee} (T_i^\vee)^{-1} = Y^{s^\vee_i\lambda^\vee}, &\textrm{if } \langle
\lambda^\vee,\alpha_i\rangle = 1,
\end{array}
\qquad\hbox{for $i= 0,1,\ldots, n$,}
\end{equation}
where the action of $W$ on $\fh_\ZZ\oplus\ZZ d$ is as in \eqref{Waction}.
\end{thm}

\subsection{Double affine Hecke algebras}

Let $R^\vee = (R^\vee)^+ \cup(-(R^\vee)^+)$ be the set of coroots and
fix parameters $c_{\beta^\vee}$, indexed by $\beta^\vee \in R^\vee +\ZZ d$, such that
for all $w\in W$ and $\beta^\vee\in R^\vee + \ZZ d$,
\begin{equation}\label{parameters}
c_{\beta^\vee} = c_{w\beta^\vee}.
\qquad\hbox{Set\quad $t_{\beta^\vee} = q^{c_{\beta^\vee}}$\quad and\quad
$t_i = t_{\alpha^\vee_i}$.}
\end{equation}
The \emph{double affine Hecke algebra} $\widetilde{H}$ is the group
algebra $\CC\widetilde{\cB}$ of the double braid group with the additional
relations
\begin{equation}\label{Heckerelation}
T_i^2 = (t_i^{1/2}-t_i^{-1/2})T_i+1, \qquad\hbox{for $i=0,1,\ldots, n$.}
\end{equation}
The double affine Hecke algebra $\widetilde H$ has bases
$$
\{ T_wX^\mu\ |\ w\in W,\  \mu\in \fh^*_\ZZ\oplus\ZZ\delta\},
\qquad
\{ Y^{\lambda^\vee}T^\vee_w\ |\ w\in W^\vee,\  \lambda^\vee\in \fh_\ZZ\oplus\ZZ d\},
$$
and
$$
\{ q^kX^\mu T_w Y^{\lambda^\vee}\ |\
w\in W_0,\  \lambda^\vee\in \fh_\ZZ, \mu\in \fh_\ZZ^*, k\in \hbox{$\frac{1}{e}$}\ZZ\}
$$
(see \cite[Prop.\ 5.4 and Cor.\ 5.8]{Hai}).

In the presence of \eqref{Heckerelation} the relations
\eqref{doubleBraidTY} are equivalent to
\begin{equation}\label{Lusztigrel}
T_i^\vee Y^{\lambda^\vee} =
Y^{s_i\lambda^\vee}T_i^\vee+(t_i^{\frac12}-t_i^{-\frac12})
\frac{Y^{\lambda^\vee}-Y^{s_i\lambda^\vee}}
{1-Y^{-\alpha_i^\vee}},
\qquad\hbox{for $i=0,1,\ldots,n$.}
\end{equation}
In turn \eqref{Lusztigrel} is equivalent to
\begin{equation}\label{Yintertwiner}
\tau_i^\vee Y^{\lambda^\vee} = Y^{s_i\lambda^\vee}\tau_i^\vee,
\qquad\hbox{for $i=0,1,\ldots, n$,}
\end{equation}
where
\begin{equation}\label{taui}
\tau_i^\vee = T_i^\vee +\frac{t_i^{-\frac12}(1-t_i)}{1-Y^{-\alpha_i^\vee}}
= (T_i^\vee)^{-1} +\frac{t_i^{-\frac12}(1-t_i)Y^{-\alpha_i^\vee}}{1-Y^{-\alpha_i^\vee}}.
\end{equation}
Using that the $\tau^\vee_i$ satisfy the braid relations and that
$$g^\vee Y^{\lambda^\vee} = Y^{g^\vee\lambda^\vee}g^\vee,
\qquad\hbox{write}\qquad
\tau^\vee_w Y^{\lambda^\vee} = Y^{w\lambda^\vee}\tau^\vee_w,
\quad\hbox{for $w \in W$.}$$

Let $w\in W$ and let $w = s^\vee_{i_1}\cdots s^\vee_{i_\ell}$ be a reduced word for $w$.
For $k=1,\ldots, \ell$ let
\begin{equation}\label{betadefn}
\beta_k^\vee = s_{i_\ell}^\vee s_{i_{\ell-1}}^\vee \cdots s_{i_{k+1}}^\vee \alpha_{i_k}^\vee
\qquad\hbox{and}\qquad t_{\beta_k^\vee}=t_{i_k},
\end{equation}
so that the sequence $\beta_\ell^\vee, \beta_{\ell-1}^\vee,
\ldots, \beta_1^\vee$ is the sequence of
labels of the hyperplanes crossed by the walk
$w^{-1} = s^\vee_{i_\ell}s^\vee_{i_{\ell-1}}\cdots s^\vee_{i_1}$.  For example,
in Type $A_2$, with
$w = s_2^\vee s_0^\vee s_1^\vee s_2^\vee s_1^\vee s_0^\vee s_2^\vee s_1^\vee$
the picture is
$$\beginpicture
\setcoordinatesystem units <1.1cm,1.1cm>         
\setplotarea x from -5 to 4, y from -2 to 4  
    \plot 1.3 -1.212 2.2 0.346 /
    \plot 0.3 -1.212 2.2 2.078 /
    \plot -0.7 -1.212 2.2 3.81 /
    \plot -1.7 -1.212 1.2 3.81 /
    \plot -2.7 -1.212 0.2 3.81 /
    \plot -3.7 -1.212 -0.8 3.81 /
    \plot -3.7 0.520 -1.8 3.81 /
    \plot -3.7 2.2516 -2.8 3.81 /
    \plot -3.3 -1.212 -3.7 -0.52 /
    \plot -2.3 -1.212 -3.7 1.212 /
    \plot -1.3 -1.212 -3.7 3.044 /
    \plot -0.3 -1.212 -3.2 3.81 /
    \plot 0.7 -1.212 -2.2 3.81 /
    \plot 1.7 -1.212 -1.2 3.81 /
    \plot 2.2 -0.346 -0.2 3.81 /
    \plot 2.2 1.386 0.8 3.81 /
    \plot 2.2 3.118 1.8 3.81 /
    \plot -3.7 -0.866 2.2 -0.866 /
    \plot -3.7 0 2.2 0 /
    \plot -3.7 3.464 2.2 3.464 /
    \plot -3.7 2.598 2.2 2.598 /
    \plot -3.7 1.732 2.2 1.732 /
    \plot -3.7 0.866 2.2 0.866 /
    \arrow <5pt> [.2,.67] from -2 2.309 to -1.5 2.021   %
    \arrow <5pt> [.2,.67] from -1.5 2.021 to -1.5 1.443   %
    \arrow <5pt> [.2,.67] from -1.5 1.443 to -1 1.155   %
    \arrow <5pt> [.2,.67] from -1 1.155 to -1 0.577   %
    \arrow <5pt> [.2,.67] from -1 0.577 to -0.5 0.289   %
    \arrow <5pt> [.2,.67] from -0.5 0.289 to 0 0.577   %
    \arrow <5pt> [.2,.67] from 0 0.577 to 0.5 0.289   %
    \arrow <5pt> [.2,.67] from 0.5 0.289 to 1 0.577   %
    \put{$\bullet$} at -2 1.732
    \put{$\scriptstyle{1}$} at -2.1 2.3
    \put{$\scriptstyle{w^{-1}}$} at 1.15 0.75
    \put{$\fh^{\beta^\vee_8}$} at -3.8 -1.5
   \put{$\fh^{\beta^\vee_7}$} at -4.1 1.732
    \put{$\fh^{\beta^\vee_6}$} at -2.8 -1.5
   \put{$\fh^{\beta^\vee_5}$} at -4.1 0.866
    \put{$\fh^{\beta^\vee_4}$} at -1.8 -1.5
    \put{$\fh^{\beta^\vee_3}$} at -0.8 -1.5
    \put{$\fh^{\beta^\vee_2}$} at 0.8 -1.5
    \put{$\fh^{\beta^\vee_1}$} at 1.8 -1.5
\endpicture$$

Let $v\in W^\vee$.
An \emph{alcove walk} of type $i_1,\ldots, i_\ell$ beginning at $v$
is a sequence of steps, where a step of type $j$ is
$$
\begin{matrix}
\beginpicture
\setcoordinatesystem units <1cm,1cm>         
\setplotarea x from -1.5 to 1.5, y from -0.5 to 0.5  
\put{$\scriptstyle{zs_j}$} at 0.6 -0.25
\put{$\scriptstyle{z}$} at -0.6 -0.25
\put{$\scriptstyle{-}$}[b] at -0.4 0.25
\put{$\scriptstyle{+}$}[b] at 0.4 0.25
\plot  0 -0.4  0 0.5 /
\arrow <5pt> [.2,.67] from -0.5 0 to 0.5 0   %
\endpicture
&\beginpicture
\setcoordinatesystem units <1cm,1cm>         
\setplotarea x from -1.5 to 1.5, y from -0.5 to 0.5  
\put{$\scriptstyle{zs_j}$} at 0.6 -0.25
\put{$\scriptstyle{z}$} at -0.6 -0.25
\put{$\scriptstyle{-}$}[b] at -0.4 0.25
\put{$\scriptstyle{+}$}[b] at 0.4 0.25
\plot  0 -0.4  0 0.5 /
\arrow <5pt> [.2,.67] from 0.5 0 to -0.5 0   %
\endpicture
&
\beginpicture
\setcoordinatesystem units <1cm,1cm>         
\setplotarea x from -1.5 to 0.5, y from -0.5 to 0.5  
\put{$\scriptstyle{zs_j}$} at 0.6 -0.25
\put{$\scriptstyle{z}$} at -0.6 -0.25
\put{$\scriptstyle{-}$}[b] at -0.4 0.35
\put{$\scriptstyle{+}$}[b] at 0.4 0.35
\plot  0 -0.4  0 0.6 /
\plot 0.5 0  0.05 0 /
\arrow <5pt> [.2,.67] from 0.05 0.1 to 0.5 0.1   %
\plot 0.05 0 0.05 0.1 /
\endpicture
&
\beginpicture
\setcoordinatesystem units <1cm,1cm>         
\setplotarea x from -1.5 to 0.5, y from -0.5 to 0.5  
\put{$\scriptstyle{zs_j}$} at 0.6 -0.25
\put{$\scriptstyle{z}$} at -0.6 -0.25
\put{$\scriptstyle{-}$}[b] at -0.4 0.35
\put{$\scriptstyle{+}$}[b] at 0.4 0.35
\plot  0 -0.4  0 0.6 /
\plot -0.5 0  -0.05 0 /
\plot -0.05 0 -0.05 0.1 /
\arrow <5pt> [.2,.67] from -0.05 0.1 to -0.5 0.1   %
\endpicture
\\
\hbox{positive $j$-crossing} &\hbox{negative $j$-crossing}
&\hbox{positive $j$-fold} &\hbox{negative $j$-fold}
\end{matrix}
$$
Let $\cB(v,\vec w)$ be the set of alcove walks of type $\vec w=(i_1,\ldots, i_\ell)$ beginning
at $v$.
For a walk $p\in \cB(v, \vec w)$ let
\begin{equation}\label{folds}
\begin{array}{rl}
f^+(p) &= \{k\ |\ \hbox{the $k$th step of $p$ is a positive fold}\},
\qquad \\
f^-(p) &= \{k\ |\ \hbox{the $k$th step of $p$ is a negative fold}\},
\end{array}
\end{equation}
and
\begin{equation}\label{endpoint}
\mathrm{end}(p) = \hbox{endpoint of $p$} \quad
\hbox{(an element of $W$).}
\end{equation}

\begin{thm} \label{mainthm}
Let $v, w\in W$, let $w = s_{i_1}^\vee\cdots s_{i_\ell}^\vee$ be a reduced
word for $w$ and let $\beta_\ell^\vee,\ldots, \beta_1^\vee$ be as defined in \eqref{betadefn}.
Then, in $\widetilde{H}$,
$$X^v \tau^\vee_w
= \sum_{p\in \cB(v,\vec w)} X^{\mathrm{end}(p)}
\left(
\prod_{k\in f^+(p)} \frac{t_{\beta_k^\vee}^{-1/2}(1-t_{\beta_k^\vee})
}{1-Y^{-\beta_k^\vee}}
\right)\left(
\prod_{k\in f^-(p)} \frac{t_{\beta_k^\vee}^{-1/2}(1-t_{\beta_k^\vee})Y^{-\beta_k^\vee}}
{1-Y^{-\beta_k^\vee}}
\right),
$$
where the sum is over all alcove walks of type $\vec w = (i_1,\ldots, i_\ell)$ beginning at $v$.
\end{thm}
\begin{proof}
The proof is by induction on the length of $w$, the base case being the formulas in
\eqref{taui}.  To do the induction step let $p\in \cB(v, \vec w)$,
$$
F^+(p) =
\left(
\prod_{k\in f^+(p)} \frac{t_{\beta_k^\vee}^{-1/2}(1-t_{\beta_k^\vee})}
{1-Y^{-\beta_k^\vee}}
\right),
\qquad
F^-(p) = \left(
\prod_{k\in f^-(p)} \frac{t_{\beta_k^\vee}^{-1/2}(1-t_{\beta_k^\vee})Y^{-\beta_k^\vee}}
{1-Y^{-\beta_k^\vee}}
\right)
$$
and let
$$\hbox{$p_1,p_2\in \cB(v,\vec ws_j)$ be the two extensions of $p$ by a step of type $j$}$$
(by a crossing and a fold, respectively).  Let $z = \mathrm{end}(p)$.
By induction, a term in $X^v\tau_w\tau_j$ is
\begin{align*}
&X^zF^+(p)F^-(p)\tau_j
=
X^z\tau_j\big(s_jF^+(p)\big)\big(s_jF^-(p)\big) \\
&=\begin{cases}
X^z
\displaystyle{\left(T_j^\vee +\frac{t_j^{-1/2}(1-t_j)}{1-Y^{-\alpha_j^\vee}}\right)}
\big(s_jF^+(p)\big)\big(s_jF^-(p)\big),
&\hbox{if $X^{zs_j}=X^zT_j^\vee$,} \\
X^z
\displaystyle{\left(
(T_j^\vee)^{-1} +\frac{t_j^{-1/2}(1-t_j)Y^{-\alpha_j^\vee}}{1-Y^{-\alpha_j^\vee}}\right)}
\big(s_jF^+(p)\big)\big(s_jF^-(p)\big),
&\hbox{if $X^{zs_j}=X^{z}(T_j^\vee)^{-1}$,}
\end{cases}
\\
&=X^{\mathrm{end}(p_1)}F^+(p_1)F^-(p_1)
+ X^{\mathrm{end}(p_2)}F^+(p_2)F^-(p_2).
\end{align*}
The last step of $p_2$ is
$$
\beginpicture
\setcoordinatesystem units <1cm,1cm>         
\setplotarea x from -1.5 to 0.5, y from -0.5 to 0.5  
\put{$\scriptstyle{-}$}[b] at -0.4 0.35
\put{$\scriptstyle{+}$}[b] at 0.4 0.35
\put{$\scriptstyle{z}$}[b] at 0.6 -0.2
\put{$\scriptstyle{zs_j}$}[b] at -0.6 -0.2
\plot  0 -0.4  0 0.6 /
\plot 0.5 0  0.05 0 /
\arrow <5pt> [.2,.67] from 0.05 0.1 to 0.5 0.1   %
\plot 0.05 0 0.05 0.1 /
\endpicture
\quad
\hbox{if $X^{zs_j} = X^zT_j^\vee$,}
\qquad\hbox{and}\quad
\beginpicture
\setcoordinatesystem units <1cm,1cm>         
\setplotarea x from -1.5 to 0.5, y from -0.5 to 0.5  
\put{$\scriptstyle{-}$}[b] at -0.4 0.35
\put{$\scriptstyle{+}$}[b] at 0.4 0.35
\put{$\scriptstyle{zs_j}$}[b] at 0.6 -0.2
\put{$\scriptstyle{z}$}[b] at -0.6 -0.2
\plot  0 -0.4  0 0.6 /
\plot -0.5 0  -0.05 0 /
\plot -0.05 0 -0.05 0.1 /
\arrow <5pt> [.2,.67] from -0.05 0.1 to -0.5 0.1   %
\endpicture
\quad
\hbox{if $X^{zs_j} = X^z (T_j^\vee)^{-1}$.}
$$
\end{proof}

\section{Macdonald polynomials}

In this section we use Theorem \ref{mainthm} to give expansions of the
nonsymmetric Macdonald polynomials $E_\mu$ (Theorem \ref{Eexpansion})
and the symmetric Macdonald polynomials $P_\mu$ (Theorem \ref{Pexpansion}).

Let $\widetilde{H}$ be the double affine Hecke algebra (defined in \eqref{Heckerelation})
and let
$H$ be the subalgebra of $\widetilde H$ generated by
$T_0,\ldots, T_n$ and $\Omega$. The \emph{polynomial representation}
of $\widetilde H$ is
\begin{equation}\label{CXdefn}
\CC[X] = \Ind_{H}^{\widetilde H}(\mathbf{1})
= \hbox{$\CC$-span}\{ q^kX^\mu\mathbf{1}\ |\ k\in \hbox{$\frac{1}{e}$} \ZZ, \mu\in \fh_\ZZ^*\}
\end{equation}
with
\begin{equation}
T_i\mathbf{1} = t_i^{1/2}\mathbf{1}
\qquad\hbox{and}\qquad g\mathbf{1} = \mathbf{1},\quad\hbox{for $g\in \Omega$.}
\end{equation}

The monomials
$X^\mu\mathbf{1}$,  $\mu\in \fh_\ZZ^*$,
form  a $\CC[q^{\pm 1/e}]$-basis of $\CC[X]$.
Another favourite $\CC[q^{\pm 1/e}]$-basis of $\CC[X]$ is
the basis of \emph{nonsymmetric Macdonald polynomials}
\begin{equation}
\{E_\mu\ |\ \mu\in \fh_\ZZ^*\},
\qquad\hbox{where\quad $E_\mu = \tau^\vee_{X^\mu m}\mathbf{1}$}
\end{equation}
with $X^\mu m$ the minimal length element in the coset $X^\mu W_0$.
Note that $\tau_w^\vee\mathbf{1}=0$ for $w\in W_0$ since
$\tau^\vee_i\mathbf{1}=0$ for $i=1,2,\ldots, n$.

If $\fh_\ZZ^+$ is the set of dominant integral coweights
(analogous to $(\fh^*_\ZZ)^+$ defined in \eqref{domint}),
$\lambda^\vee\in \fh_\ZZ^+$ and
$Y^{\lambda^\vee} = s_{i_1}\cdots s_{i_\ell}$ is a reduced word,
then
$$
Y^{\lambda^\vee} \mathbf{1}
= T_{i_1}\cdots T_{i_\ell} \mathbf{1}
= t_{i_1}^{\frac12}\cdots t_{i_\ell}^{\frac12}\mathbf{1}
= q^{\frac12(c_{i_1}+\cdots+c_{i_\ell})}\mathbf{1}
= q^{\frac12\sum_{\alpha\in R^+} c_\alpha\langle
\lambda^\vee,\alpha\rangle}\mathbf{1}
=q^{\langle \lambda^\vee,\rho_c\rangle}\mathbf{1},
$$
since $\langle \lambda^\vee,\alpha\rangle$ is the number of
hyperplanes parallel to $\fh^\alpha$ which are between
$Y^{\lambda^\vee}$ and $1$.  If $\lambda^\vee\in \fh_\ZZ$ then
$\lambda^\vee = \mu^\vee-\nu^\vee$ for some
$\mu^\vee$, $\nu^\vee \in \fh_\ZZ^+$ and so, for all
$\lambda^\vee \in \fh_\ZZ$,
\begin{equation}\label{Y1}
Y^{\lambda^\vee} \mathbf{1} =
q^{\langle \lambda^\vee,\rho_c\rangle}\mathbf{1},
\qquad\hbox{where}\quad \rho_c = \hbox{$\frac12$}\sum_{\alpha\in
R^+} c_\alpha \alpha.
\end{equation}
More generally, if $X^\mu m$ is the minimal length element of the coset
$X^\mu W_0$
then
\begin{align*}
Y^{\lambda^\vee}E_\mu
&= Y^{\lambda^\vee}\tau_{X^\mu m}\mathbf{1}
=\tau_{X^\mu m}Y^{m^{-1} X^{-\mu}\lambda^\vee}\mathbf{1}
=\tau_{X^\mu m} Y^{m^{-1}(\lambda^\vee + \langle \lambda^\vee,\mu\rangle d)}\mathbf{1} \\
&=\tau_{X^\mu m} Y^{m^{-1}\lambda^\vee}q^{-\langle \lambda^\vee,\mu\rangle}\mathbf{1}
=q^{\langle m^{-1}\lambda^\vee,\rho_c\rangle-\langle \lambda^\vee,\mu\rangle}
\tau_{X^\mu m}\mathbf{1}
=q^{\langle \lambda^\vee,m\rho_c-\mu\rangle}
E_\mu \\
&=q^{\langle \lambda^\vee,X^{-\mu}m\cdot \rho_c \rangle}
E_\mu,
\end{align*}
where, in the last line, the action of $W$ on $\fh_\ZZ^*$ is as in \eqref{translaction}.  Thus the
$E_\mu$ are eigenvectors for the action of the $Y^{\lambda^\vee}$ on the polynomial
representation $\CC[X]$.

Retain the notation of (\ref{folds}-\ref{endpoint}) so that if
$w = s^\vee_{i_1}\cdots s^\vee_{i_\ell}$ is a reduced word then
$\cB(v,\vec w)$ denotes the
set of alcove walks of type $\vec w = (i_1,\ldots, i_\ell)$ beginning at $v$.
For $p\in \cB(v,\vec w)$ define the \emph{weight $\wt(p)$} and the
\emph{final direction $\varphi(p)$} of $p$ by
\begin{equation}\label{findir}
X^{\mathrm{end}(p)} = X^{\wt(p)} T^\vee_{\varphi(p)},\qquad
\hbox{with $\wt(p)\in \fh_\ZZ^*$ and $\varphi(p)\in W_0$.}
\end{equation}
In other words, $\wt(p)$ is the ``hexagon where $p$ ends''.  For $w\in W$ define
\begin{equation}\label{tlength}
t^{1/2}_w = t^{1/2}_{i_1}\ldots t^{1/2}_{i_\ell},
\qquad\hbox{if $w = s^\vee_{i_1}\cdots s^\vee_{i_\ell}$ is a reduced word.}
\end{equation}
If $\beta_k^\vee = s_{i_\ell}^\vee \cdots s_{i_{k+1}}^\vee\alpha_{i_k}$ are as defined in
\eqref{betadefn} then, by \eqref{Y1},
$$Y^{-\beta_k^\vee}\mathbf{1} = Y^{-(-\gamma^\vee+jd)}\mathbf{1}
= q^jq^{\langle \gamma^\vee, \rho_c\rangle}\mathbf{1},
\qquad\hbox{if $\beta_k^\vee = -\gamma^\vee+jd$}
$$
with $\gamma^\vee\in R^\vee$, $j\in \ZZ$.  By \eqref{parameters} and the definition of
$\rho_c$ in \eqref{Y1},
the constant $q^{\langle \gamma^\vee, \rho_c\rangle}$ is a monomial in the
symbols $t_i^{1/2}$.  To simplify
the notation for these constants write $q^jq^{\langle \gamma^\vee,\rho_c\rangle}
=q^{\langle-\beta_k^\vee,\rho_c\rangle}$ so that
\begin{equation}\label{constants}
Y^{-\beta_k^\vee}\mathbf{1} = q^{\langle -\beta_k^\vee, \rho_c\rangle}\mathbf{1}.
\end{equation}

\begin{thm}\label{Eexpansion}
Let $\mu\in \fh_\ZZ^*$ and let $w = X^\mu m$ be the minimal length element
in the coset $X^\mu W_0$.  Fix a reduced word $\vec w = s_{i_1}^\vee\cdots s_{i_\ell}^\vee$
for $w$ and let $\beta_\ell^\vee,\ldots, \beta_1^\vee$ be as defined in \eqref{betadefn}.
With notations as in (\ref{findir}-\ref{constants}) the nonsymmetric Macdonald polynomial
$$
E_\mu
= \sum_{p\in \cB(\vec \mu)} X^{\wt(p)} t^{\frac12}_{\varphi(p)}
\left(
\prod_{k\in f^+(p)} \frac{t_{\beta_k^\vee}^{-\frac12}(1-t_{\beta_k^\vee}) }
{1-q^{\langle -\beta_k^\vee,\rho_c\rangle}}
\right)\left(
\prod_{k\in f^-(p)} \frac{t_{\beta_k^\vee}^{-\frac12}(1-t_{\beta_k^\vee})
q^{\langle-\beta_k^\vee,\rho_c\rangle}}
{1-q^{\langle-\beta_k^\vee,\rho_c\rangle}}
\right),
$$
where the sum is over the set $\cB(\vec \mu)=\cB(1, \vec w)$ of alcove walks of type
$i_1,\ldots, i_\ell$ beginning at $1$.
\end{thm}
\begin{proof}
Since $E_\mu = \tau^\vee_{X^\mu m}\mathbf{1}$,
$$X^{\mathrm{end}(p)}\mathbf{1} = X^{\wt(p)}T_{\varphi(p)}^\vee\mathbf{1}
=X^{\wt(p)}t^{\frac12}_{\varphi(p)}\mathbf{1}
\qquad\hbox{and}\qquad
Y^{\lambda^\vee}\mathbf{1} = q^{\langle \lambda^\vee, \rho_c\rangle}\mathbf{1},
$$
applying the formula for $\tau_{X^\mu m}$
in Theorem \ref{mainthm} to $\mathbf{1}$ gives
the formula in the statement.
\end{proof}

\begin{remark} From the expansion of $E_\mu$ in Theorem \ref{Eexpansion},
the nonsymmetric Macdonald
polynomial $E_\mu$ has top term $t^{1/2}_m X^\mu$, where $X^\mu m$ is the minimal
length representative of the coset $X^\mu W_0$.  This term is the term corresponding to the
unique alcove walk in $\cB(\vec \mu)$ with no folds.
\end{remark}

\begin{remark}  If $w=X^\mu m = s^\vee_{i_1}\cdots s_{i_\ell}^\vee$ is a reduced word
for the minimal length element of the coset $X^\mu W_0$ then
$w^{-1} = s_{i_\ell}^\vee \cdots s^\vee_{i_1}$ is a walk from $1$ to $w^{-1}$ which stays
completely in the dominant chamber.  This has the effect that the roots
$\beta_\ell^\vee, \ldots, \beta_1^\vee$ are all of the form $-\gamma^\vee+jd$ with
$\gamma^\vee\in (R^\vee)^+$ (positive coroots) and $j\in \ZZ_{> 0}$.  The \emph{height}
of a coroot $\gamma^\vee$ is
$$\mathrm{ht}(\gamma^\vee) = \langle \gamma^\vee, \rho\rangle,
\qquad\hbox{where}\ \ \rho = \hbox{$\frac12$}\sum_{\alpha\in R} \alpha.
$$
In the case that all the parameters are equal ($t_i=t=q^c$ for $i=0,\ldots, n$)
the values which appear in Theorem \ref{Eexpansion},
$$q^{\langle -\beta_k^\vee, \rho_c\rangle}
= q^{\langle \gamma^\vee-jd,\rho_c\rangle}
= q^j t^{\mathrm{ht}(\gamma^\vee)},
\qquad\hbox{
 have positive exponents (in $\ZZ_{>0}$).}
 $$
\end{remark}

The set of \emph{dominant integral weights} is
\begin{equation}\label{domint}
(\fh_\ZZ^*)^+
= \{ \mu\in \fh_\ZZ^*\ |\
\hbox{$\langle \mu,\alpha_i^\vee\rangle \ge 0$ for
$i=1,\ldots, n$} \}.
\end{equation}
Recall the notation for $t_w^{1/2}$ from \eqref{tlength}.
For $\mu\in (\fh_\ZZ^*)^+$, the \emph{symmetric Macdonald polynomial}
(see \cite[Remarks after (6.8)]{Mac3}) is
\begin{equation}
P_\mu = \mathbf{1}_0 E_\mu
\qquad\hbox{where}\quad
\mathbf{1}_0 = \sum_{w\in W_0} t_{w_0w}^{-\frac12}T_w,
\end{equation}
so that
$T_i \mathbf{1}_0 = t_i^{1/2}\mathbf{1_0}$ for $i=1,2\ldots n$,
and $\mathbf{1}_0$ has top term $T_{w_0}$ with coefficient $1$.  The symmetric
Macdonald polynomials are $W_0$-symmetric polynomials in $X^\mu$ which are
eigenvectors for the action of $W_0$-symmetric polynomials in the $Y^{\lambda^\vee}$.

\begin{thm}\label{Pexpansion}
Let $\mu\in (\fh_\ZZ^*)^+$ and let $X^\mu m = s^\vee_{i_1}\cdots s^\vee_{i_\ell}$
be a reduced word for the minimal length element $X^\mu m$ in the coset $X^\mu W_0$.
Let $\beta_\ell^\vee,\ldots, \beta_1^\vee$ be as defined in \eqref{betadefn} and let
$$\cP(\vec\mu) = \bigcup_{v\in W_0} \cB(v,\vec w)$$
be the set of alcove walks of type $\vec{w} = (i_1,\ldots, i_\ell)$ beginning
at an element $v\in W_0$.
Then the symmetric Macdonald polynomial
$$P_\mu = \sum_{p\in \cP(\vec \mu)}
X^{\wt(p)} t^{\frac12}_{\varphi(p)}t^{-\frac12}_{w_0w}
\left(
\prod_{k\in f^+(p)} \frac{t_{\beta_k^\vee}^{-\frac12}(1-t_{\beta_k^\vee})}
{1-q^{\langle -\beta_k^\vee,\rho_c\rangle}}
\right)\left(
\prod_{k\in f^-(p)} \frac{t_{\beta_k^\vee}^{-\frac12}(1-t_{\beta_k^\vee})
q^{\langle-\beta_k^\vee,\rho_c\rangle}}
{1-q^{\langle-\beta_k^\vee,\rho_c\rangle}}
\right).
$$
\end{thm}
\begin{proof}
The expression
$$\mathbf{1}_0 = \sum_{w\in W_0} t_{w_0w}^{-\frac12}X^w,
\qquad\hbox{gives}\qquad
P_\mu = \mathbf{1}_0 E_{\mu} =
\sum_{w\in W_0} t^{-\frac12}_{w_0w}X^w
\tau^\vee_{X^\mu m} \mathbf{1},$$
which is computed by the same method as in Theorem \ref{mainthm}
and Theorem \ref{Eexpansion}.
\end{proof}

\begin{remark}
The \emph{Hall-Littlewood polynomials} or \emph{Macdonald spherical functions}
are $P_\mu(0,t)$ and the \emph{Schur functions} or \emph{Weyl characters}  are
$s_\mu = P_\mu(0,0)$.   In the first case the formula
in Theorem \ref{Pexpansion} reduces to the formula for the Macdonald spherical functions
in terms of positively folded alcove walks as given in \cite[Thm.\ 1.1]{Sc}
(see also \cite[Thm.\ 4.2(a)]{Ra}).
In the case $q=t=0$, the formula in Theorem \ref{Pexpansion} reduces
to the formula for the Weyl characters in terms of maximal dimensional positively
folded alcove walks (\emph{the Littelmann path model}) as given in
\cite[Cor.\ 1 p.\ 62]{GL}.
\end{remark}

\end{section}

\section{Examples}

\subsection{Type $A_1$}
The Weyl group $W_0=\langle s_1\ |\ s_1^2=1\rangle$ has order two
and acts on the lattices
\begin{equation}
\fh_\ZZ = \ZZ \omega^\vee
\quad\hbox{and}\quad \fh_\ZZ^* = \ZZ\omega \qquad\hbox{by}\qquad
s_1\omega^\vee = -\omega^\vee \quad\hbox{and}\quad s_1\omega = -
\omega,
\end{equation}
and
\begin{equation}
\varphi^\vee = \alpha^\vee = 2\omega^\vee,
\qquad \varphi = \alpha = 2\omega, \quad\hbox{and}\quad \langle
\omega^\vee,\alpha\rangle = 1.
\end{equation}
$$\beginpicture
\setcoordinatesystem units <1.5cm,1cm>         
\setplotarea x from -4 to 5, y from -1 to 3  
\put{$\bullet$} at 0 0 \put{$\fh^{\alpha^\vee-2d}$}[t] at 2 -0.7
\put{$\fh^{\alpha^\vee}$}[t] at 0 -0.7 \put{$\fh^{\alpha^\vee+2d}$}[t] at -2
-0.7 \put{$\fh^{\alpha^\vee-d}$}[t] at 1 -0.7
\put{$\fh^{\alpha^\vee-3d}$}[t] at 3 -0.7 \put{$\fh^{\alpha^\vee+d}$}[t]
at -1 -0.7
\plot -4 0  5 0 / \plot  -3 0.3  -3 -0.3 / \plot  -2 0.3  -2 -0.3 /
\plot  -1 0.3  -1 -0.3 / \plot  0 0.3  0 -0.3 / \plot  1 0.3  1 -0.3
/ \plot  2 0.3  2 -0.3 / \plot  3 0.3  3 -0.3 / \plot  4 0.3  4 -0.3
/
\small{
\put{$1$}[b] at 0.5 0.3 \put{$s_0^\vee$}[b] at 1.5 0.3 \put{$s_0^\vee s_1^\vee$}[b]
at 2.5 0.3 \put{$s_0^\vee s_1^\vee s_0^\vee$}[b] at 3.5 0.3 \put{$s_1^\vee$}[b] at -0.5
0.3 \put{$s_1^\vee s_0^\vee$}[b] at -1.5 0.3 \put{$s_1^\vee s_0^\vee s_1^\vee$}[b] at -2.5 0.3
\put{$1$}[t] at 0.5 -0.3 \put{$X^{\alpha}s_1$}[t] at 1.5 -0.3
\put{$X^{\alpha}$}[t] at 2.5 -0.3
\put{$X^{2\alpha}s_1$}[t] at 3.5 -0.3 \put{$s_1$}[t] at -0.5
-0.3 \put{$X^{-\alpha}$}[t] at -1.5 -0.3
\put{$X^{-\alpha}s_1$}[t] at -2.5 -0.3
\plot -4 2 5 2 / \plot  -3 2.3  -3 1.7 / \plot  -2 2.3  -2 1.7 /
\plot  -1 2.3  -1 1.7 / \plot  0 2.3  0 1.7 / \plot  1 2.3  1 1.7 /
\plot  2 2.3  2 1.7 / \plot  3 2.3  3 1.7 / \plot  4 2.3  4 1.7 /
\put{$g^\vee$}[b] at 0.5 2.3 \put{$g^\vee s_1^\vee$}[b] at 1.5 2.3
\put{$g^\vee s_1^\vee s_0^\vee$}[b] at 2.5 2.3 \put{$g^\vee s_1^\vee s_0^\vee s_1^\vee$}[b] at 3.5 2.3
\put{$g^\vee s_0^\vee$}[b] at -0.5 2.3 \put{$g^\vee s_0^\vee s_1^\vee$}[b] at -1.5 2.3
\put{$g^\vee s_0^\vee s_1^\vee s_0^\vee$}[b] at -2.5 2.3
\put{$X^{\omega}s_1$}[t] at 0.5 1.7 \put{$X^{\omega}$}[t]
at 1.5 1.7 \put{$X^{3\omega}s_1$}[t] at 2.5 1.7
\put{$X^{3\omega}$}[t] at 3.5 1.7 \put{$X^{-\omega}$}[t]
at -0.5 1.7 \put{$X^{-\omega}s_1$}[t] at -1.5 1.7
\put{$X^{-3\omega}$}[t] at -2.5 1.7
} 
\endpicture
$$

The double affine braid group $\widetilde \cB$ is generated by $T_0, T_1,
g, X^\omega$, and $q^{1/2}$, with relations
\begin{equation} \label{A1dblaffbdgp}
\begin{array}{c}
T_0=gT_1g^{-1},\qquad g^2=1,\qquad q = X^\delta, \\ \\
gX^\omega = q^{1/2}X^{-\omega}g, \qquad T_1X^\omega T_1=X^{-\omega},
\quad\textrm{and}\quad T_0X^{-\omega}T_0=q^{-1}X^\omega.
\end{array}
\end{equation}
In the double affine braid group
\begin{equation}\label{A1altgen}
g=Y^{\omega^\vee}T_1^{-1},
\qquad
T_0 = Y^{\varphi^\vee} T_1^{-1},
\qquad
g^\vee =X^\omega T_1,
\qquad
(T_0^\vee)^{-1}=X^{\varphi}T^\vee_1.
\end{equation}
At this point, the following Proposition, which is the Type $A_1$ case of Theorem \ref{duality},
is easily proved by direct computation.

\begin{prop}\label{A1duality}
(Duality). Let $Y^d = q^{-1}$. The double affine braid group
$\widetilde{\cB}$ is generated by $T_0^\vee, T_1^\vee, g^\vee,
Y^{\omega^\vee}$ and $q^{1/2}$ with relations
$$Y^d=q^{-1},
\qquad
(g^\vee)^2 =1,\qquad
T_0^\vee = g^\vee T_1^\vee(g^\vee)^{-1},$$
$$g^\vee Y^{\omega^\vee} =q^{-1/2}Y^{-\omega^\vee}g^\vee,
\quad T_1^{-1}Y^{\omega^\vee}T_1^{-1} = Y^{-\omega^\vee},
\quad\hbox{and}\quad (T_0^\vee)^{-1} Y^{-\omega^\vee}(T_0^\vee)^{-1} = qY^{\omega^\vee}.
$$
\end{prop}

The double affine Hecke algebra $\widetilde H$ is
$\CC\widetilde{\cB}$ with the additional relations
\begin{equation}\label{A1Hecke}
T_i^2 = (t^{1/2}-t^{-1/2})T_i+1,\ \ \hbox{for $i=0,1$,}
\qquad\hbox{and}\qquad
t_0=t_1=t=q^c.
\end{equation}
Using \eqref{A1Hecke}, the relations in Proposition
\eqref{A1duality} give
$$
g^\vee Y^{\omega^\vee} =q^{-1/2}Y^{-\omega^\vee}g^\vee, \qquad
T_1Y^{\omega^\vee} =Y^{-\omega^\vee}T_1 +
(t^{1/2}-t^{-1/2})\frac{Y^{\omega^\vee}-Y^{-\omega^\vee}}{1-Y^{-\alpha^\vee}},
\quad\hbox{and}
$$
$$
T^\vee_0 Y^{\omega^\vee} =q^{-1}Y^{-\omega^\vee}T^\vee_0
+(t^{1/2}-t^{-1/2})
\left(\frac{Y^{\omega^\vee}-q^{-1}Y^{-\omega^\vee}}
{1-qY^{\alpha^\vee}}\right).$$
With $Y^{\alpha_0} = qY^\alpha$ and $Y^{\alpha_1}=Y^\alpha$, then
$$\tau^\vee_g = g^\vee, \qquad\hbox{and}\qquad
\tau^\vee_i = T_i^\vee -
(t^{1/2}-t^{-1/2})\left(\frac{1}{1-Y^{-\alpha^\vee_i}}\right),
\qquad\hbox{for $i=0,1$}.$$
%

To illustrate Theorem \ref{mainthm}, note that $X^{-2\omega}=s_1^\vee s^\vee_0$ is a
reduced word and
\begin{align*}
\tau_1^\vee &\tau_0^\vee
= \left( T_1^\vee + \frac{t^{-1/2}(1-t)}{1-Y^{-\alpha_1^\vee}}\right)\tau_0^\vee \\
&= T_1^\vee T_0^\vee + T_1^\vee \frac{t^{-1/2}(1-t)}{1-Y^{-\alpha_0^\vee}} 
+(T_0^\vee)^{-1}\frac{t^{-1/2}(1-t)}{1-Y^{-s_0\alpha_1^\vee}}
+\left(\frac{t^{-1/2}(1-t)}{1-Y^{-s_0\alpha_1^\vee}}\right)
\left(\frac{t^{-1/2}(1-t)Y^{-\alpha_0^\vee}}{1-Y^{-\alpha_0^\vee}}\right) \\
&= X^{-2\omega} + T_1^\vee \frac{t^{-1/2}(1-t)}{1-Y^{-\alpha_0^\vee}} 
+X^{2\omega}T_1^\vee \frac{t^{-1/2}(1-t)}{1-Y^{-s_0\alpha_1^\vee}}
+\left(\frac{t^{-1/2}(1-t)}{1-Y^{-s_0\alpha_1^\vee}}\right)
\left(\frac{t^{-1/2}(1-t)Y^{-\alpha_0^\vee}}{1-Y^{-\alpha_0^\vee}}\right).
\end{align*}
The corresponding
paths in $\cB(1,\overrightarrow{-2\omega})=\cB(\overrightarrow{-2\omega})$ are
$$
\begin{matrix}
\beginpicture
\setcoordinatesystem units <.7cm,.7cm>         
\setplotarea x from -2.5 to 1.5, y from -0.5 to 0.6  
    \plot -2.2 0 1.2 0 /
    \plot  -2 -0.2  -2 0.6 /
    \plot  -1 -0.2  -1 0.6 /
    \plot  0 -0.2  0 0.6 /
    \plot  1 -0.2  1 0.6 /
    \put{$\bullet$} at 0 0
    \arrow <5pt> [.2,.67] from 0.5 0.2 to -0.5 0.2
    \arrow <5pt> [.2,.67] from -0.5 0.2 to -1.5 0.2
\endpicture
&
\beginpicture
\setcoordinatesystem units <.7cm,.7cm>         
\setplotarea x from -2.5 to 1.5, y from -0.5 to 0.6  
    \plot -2.2 0 1.2 0 /
    \plot  -2 -0.2  -2 0.6 /
    \plot  -1 -0.2  -1 0.6 /
    \plot  0 -0.2  0 0.6 /
    \plot  1 -0.2  1 0.6 /
    \put{$\bullet$} at 0 0
    \arrow <5pt> [.2,.67] from 0.5 0.2 to -0.5 0.2
    \plot -0.5 0.2 -1 0.2 /
    \plot -0.95 0.2 -0.95 0.35 /
    \arrow <5pt> [.2,.67] from -1 0.35 to -0.5 0.35
\endpicture
&
\beginpicture
\setcoordinatesystem units <.7cm,.7cm>         
\setplotarea x from -2.5 to 1.5, y from -0.5 to 0.6  
    \plot -2.2 0 1.2 0 /
    \plot  -2 -0.2  -2 0.6 /
    \plot  -1 -0.2  -1 0.6 /
    \plot  0 -0.2  0 0.6 /
    \plot  1 -0.2  1 0.6 /
    \put{$\bullet$} at 0 0
    \plot 0.5 0.2 0 0.2 /
    \plot 0.05 0.2 0.05 0.35 /
    \arrow <5pt> [.2,.67] from 0 0.35 to 0.5 0.35
    \arrow <5pt> [.2,.67] from 0.5 0.35 to 1.5 0.35
\endpicture
&
\beginpicture
\setcoordinatesystem units <.7cm,.7cm>         
\setplotarea x from -2.5 to 1.5, y from -0.5 to 0.6  
    \plot -2.2 0 1.2 0 /
    \plot  -2 -0.2  -2 0.6 /
    \plot  -1 -0.2  -1 0.6 /
    \plot  0 -0.2  0 0.6 /
    \plot  1 -0.2  1 0.6 /
    \put{$\bullet$} at 0 0
    \plot 0.5 0.2 0 0.2 /
    \plot 0.05 0.2 0.05 0.35 /
    \arrow <5pt> [.2,.67] from 0 0.35 to 0.5 0.35
    \plot 0.5 0.35 1 0.35 /
    \plot 0.95 0.35 0.95 0.5 /
    \arrow <5pt> [.2,.67] from 1 0.5 to 0.5 0.5
\endpicture
\\ \\
X^{-2\omega}
&\displaystyle{ T_1^\vee \frac{t^{-1/2}(1-t)}{1-Y^{-\alpha_0^\vee}} }
&\displaystyle{ X^{2\omega} T_1^\vee \frac{t^{-1/2}(1-t)}{1-Y^{-s_0\alpha_1^\vee}} }
&\qquad\displaystyle{ \frac{t^{-1/2}(1-t)}{1-Y^{-s_0\alpha_1^\vee}}
\frac{t^{-1/2}(1-t)Y^{-\alpha_0^\vee}}{1-Y^{-\alpha_0^\vee}} }
\end{matrix}
$$

The polynomial representation is defined by
$$T_i\mathbf{1} = t^{1/2}\mathbf{1},
\qquad\hbox{and}\qquad g\mathbf{1} = \mathbf{1}.
$$
In this case
\begin{equation}
\rho_c = \frac12c\alpha
\qquad\hbox{and}\qquad
W^0 = \{ X^{-\ell\omega}\ |\ \ell\in \ZZ_{\ge 0}\}
\cup \{ X^{\ell\omega}s^\vee_1\ |\ \ell\in \ZZ_{>0}\},
\end{equation}
is the set of minimal length coset representatives of $W^\vee/W_0$.

Applying the expansion of $\tau_1^\vee\tau_0^\vee$ to $\mathbf{1}$ and using
\begin{align*}
Y^{-\alpha_0}\mathbf{1} 
= qY^{\alpha^\vee}\mathbf{1} = qq^c\mathbf{1} = qt\mathbf{1},
\qquad\hbox{and}\qquad
Y^{-s_0\alpha_1^\vee}\mathbf{1} = Y^{\alpha^\vee+2d}\mathbf{1}
= q^2Y^{\alpha^\vee}\mathbf{1} = q^2t,
\end{align*}
gives
\begin{align*}
E_{-2\omega} &= \tau_1^\vee \tau_0^\vee\mathbf{1} \\
&= X^{-2\omega} + t^{1/2} \frac{t^{-1/2}(1-t)}{1-qt}
+X^{2\omega}t^{1/2} \frac{t^{-1/2}(1-t)}{1-q^2t}
+\left(\frac{t^{-1/2}(1-t)}{1-q^2t}\right)
\left(\frac{t^{-1/2}(1-t)qt}{1-qt}\right) \\
&= X^{-2\omega} + \frac{1-t}{1-qt}
+X^{2\omega} \frac{1-t}{1-q^2t}
+\left(\frac{1-t}{1-q^2t}\right)
\left(\frac{(1-t)q}{1-qt}\right)
\end{align*}

Since
$\mathbf{1}_0 = T_1^\vee+t^{-1/2}$ the symmetric Macdonald polynomial
$P_{2\omega} = \mathbf{1}_0E_{2\omega} = \mathbf{1}_0\tau_0^\vee\mathbf{1}$ is
\begin{align*}
P_{2\omega} &= \mathbf{1_0}E_{2\omega}
= (T_1^\vee + t^{-1/2})\tau_0^\vee \mathbf{1} \\
&= \left(
T_1^\vee T_0^\vee + T_1^\vee\frac{t^{-1/2}(1-t)}{1-Y^{-\alpha_0^\vee}}
+t^{-1/2}(T_0^\vee)^{-1}
+t^{-1/2}\frac{t^{-1/2}(1-t)Y^{-\alpha_0^\vee}}{1-Y^{-\alpha_0^\vee}}\right)\mathbf{1} \\
&= \left(
X^{-2\omega} + t^{1/2}\frac{t^{-1/2}(1-t)}{1-qt}
+t^{-1/2}X^{2\omega}T_1^\vee
+t^{-1/2}\frac{t^{-1/2}(1-t)qt}{1-qt}\right)\mathbf{1} \\
&= \left(
X^{-2\omega} + \frac{1-t}{1-qt}
+X^{2\omega}
+\frac{(1-t)q}{1-qt}\right)\mathbf{1}
=\left( X^{2\omega}+X^{-2\omega}+(1+q)\frac{1-t}{1-qt}\right)\mathbf{1}.
\end{align*}
The corresponding paths in $\cP(\overrightarrow{2\omega})$ are
$$
\begin{matrix}
\beginpicture
\setcoordinatesystem units <.7cm,.7cm>         
\setplotarea x from -2.5 to 2.5, y from -0.5 to 0.5  
    \plot -2.2 0 2.2 0 /
    \plot  -2 -0.2  -2 0.5 /
    \plot  -1 -0.2  -1 0.5 /
    \plot  0 -0.2  0 0.5 /
    \plot  1 -0.2  1 0.5 /
    \plot 2 -0.2 2 0.5 /
    \put{$\bullet$} at 0 0
    \arrow <5pt> [.2,.67] from -0.5 0.2 to -1.5 0.2
\endpicture
&
\beginpicture
\setcoordinatesystem units <.7cm,.7cm>         
\setplotarea x from -2.5 to 2.5, y from -0.5 to 0.5  
    \plot -2.2 0 2.2 0 /
    \plot  -2 -0.2  -2 0.5 /
    \plot  -1 -0.2  -1 0.5 /
    \plot  0 -0.2  0 0.5 /
    \plot  1 -0.2  1 0.5 /
    \plot 2 -0.2 2 0.5 /
    \put{$\bullet$} at 0 0
    \plot -0.5 0.2 -1 0.2 /
    \plot -0.95 0.2 -0.95 0.35 /
    \arrow <5pt> [.2,.67] from -1 0.35 to -0.5 0.35
\endpicture
&
\beginpicture
\setcoordinatesystem units <.7cm,.7cm>         
\setplotarea x from -2.5 to 2.5, y from -0.5 to 0.5  
    \plot -2.2 0 2.2 0 /
    \plot  -2 -0.2  -2 0.5 /
    \plot  -1 -0.2  -1 0.5 /
    \plot  0 -0.2  0 0.5 /
    \plot  1 -0.2  1 0.5 /
    \plot 2 -0.2 2 0.5 /
    \put{$\bullet$} at 0 0
    \arrow <5pt> [.2,.67] from 0.5 0.2 to 1.5 0.2
\endpicture
&
\beginpicture
\setcoordinatesystem units <.7cm,.7cm>         
\setplotarea x from -2.5 to 2.5, y from -0.5 to 0.5  
    \plot -2.2 0 2.2 0 /
    \plot  -2 -0.2  -2 0.5 /
    \plot  -1 -0.2  -1 0.5 /
    \plot  0 -0.2  0 0.5 /
    \plot  1 -0.2  1 0.5 /
    \plot 2 -0.2 2 0.5 /
    \put{$\bullet$} at 0 0
    \plot 0.5 0.2 1 0.2 /
    \plot 0.95 0.2 0.95 0.35 /
    \arrow <5pt> [.2,.67] from 1 0.35 to 0.5 0.35
\endpicture
\\
\displaystyle{X^{-2\omega}} & \displaystyle{\frac{1-t}{1-tq}}
& \displaystyle{X^{2\omega}}
& \displaystyle{q\frac{(1-t)}{1-tq}}
\end{matrix}
$$

\subsection{Type $A_2$}
The Weyl group $W_0=\langle s_1, s_2 \mid s_1^2=s_2^2=1,
s_1s_2s_1=s_2s_1s_2\rangle$ acts on the lattices
\begin{equation}
\fh_\ZZ = \ZZ\omega_1^\vee+\ZZ\omega_2^\vee
\qquad\hbox{and}\qquad \fh_\ZZ^* = \ZZ\omega_1+\ZZ\omega_2,
\end{equation}
where $s_1$ and $s_2$ are the reflections in the hyperplanes
determined by
\begin{equation}
\alpha_1^\vee = 2\omega_1^\vee - \omega_2^\vee, \quad \alpha_2^\vee
= -\omega_1^\vee + 2\omega_2^\vee, \quad \alpha_1 = 2\omega_1 -
\omega_2, \quad\textrm{and}\quad \alpha_2 = -\omega_1 + 2\omega_2,
\end{equation}
with $\langle\omega_i^\vee,\alpha_j\rangle = \delta_{ij}$, and
$\langle\omega_i,\alpha_j^\vee\rangle = \delta_{ij}.$
In this case,
\begin{equation}
\varphi^\vee = \alpha_1^\vee+\alpha_2^\vee,
\qquad\hbox{and}\qquad
\varphi = \alpha_1+\alpha_2.
\end{equation}

The double affine braid group $\widetilde{\cB}$ is generated by $T_0, T_1,
T_2, g, X^{\omega_1}, X^{\omega_2}$, and $q^{1/3}$, with relations
\begin{equation}
\begin{array}{c}
T_iT_jT_i = T_jT_iT_j, \qquad\hbox{for $i\neq j$,}\\ \\
X^{\mu}X^{\lambda} = X^{\mu+\lambda} = X^{\lambda}X^{\mu},
\qquad  \hbox{for $\mu,\lambda \in \fh^*_\ZZ$,}\\ \\
T_1 X^{\omega_2} =X^{\omega_2}T_1, \quad T_2 X^{\omega_1}
=X^{\omega_1}T_2, \quad T_1X^{\omega_1}T_1 = X^{-\omega_1+\omega_2},
\quad T_2X^{\omega_2}T_2 = X^{\omega_1-\omega_2},\\  \\
g^3=1, \qquad gX^{\omega_1}= q^{1/3}X^{-\omega_1+\omega_2}g,
\quad gX^{\omega_2}= q^{2/3} X^{-\omega_1}g,\\ \\
gT_0g^{-1} = T_1, \qquad gT_1g^{-1} = T_2, \qquad gT_2g^{-1} = T_0.
\end{array}
\end{equation}
The formula \eqref{gT0} gives
\begin{equation}
g = Y^{\omega_1^\vee}T_1^{-1}T_2^{-1},
\qquad
g^2 = Y^{\omega_2^\vee}T_2^{-1}T_1^{-1},
\qquad
T_0 = Y^{\varphi^\vee}T_1^{-1}T_2^{-1}T_1^{-1},
\end{equation}
\begin{equation}
g^\vee = X^{\omega_1} T^\vee_1T^\vee_2,
\qquad
(g^\vee)^2 = X^{\omega_2}T_2^\vee T_1^\vee,
\qquad
(T_0^\vee)^{-1} = X^\varphi T_1^\vee T_2^\vee T_1^\vee.
\end{equation}
At this point, the following Proposition, which is the Type $A_2$ case of Theorem \ref{duality},
is easily proved by direct computation.

\begin{prop}(Duality). Let $Y^d = q^{-1}$.
The double affine braid group $\widetilde{\cB}$ is generated by $T_0^\vee$,
$T_1^\vee$, $T_2^\vee$, $g^\vee$, $Y^{\omega_1^\vee}$,
$Y^{\omega_2^\vee}$ and $q^{1/3}$, with relations
$$(T_1^\vee)^{-1} Y^{\omega_1^\vee} (T_1^\vee)^{-1} = Y^{-\omega_1^\vee+\omega_2^\vee},
\quad (T_2^\vee)^{-1} Y^{\omega_2^\vee} (T_2^\vee)^{-1} = Y^{\omega_1^\vee-\omega_2^\vee},$$
$$(T_1^\vee)^{-1} Y^{\omega_2^\vee} =Y^{\omega_2^\vee}(T_1^\vee)^{-1},
\quad (T_2^\vee)^{-1} Y^{\omega_1^\vee} =Y^{\omega_1^\vee}(T_2^\vee)^{-1},$$
$$(g^\vee)^3=1, \quad g^\vee Y^{\omega_1^\vee} = q^{-1/3}Y^{-\omega_1^\vee+\omega_2^\vee}g^\vee,
\quad g^\vee Y^{\omega_2^\vee} = q^{-2/3}Y^{-\omega_1^\vee}g^\vee,$$
$$g^\vee T_0^\vee (g^\vee)^{-1}= T_1^\vee, \qquad
g^\vee T_1^\vee (g^\vee)^{-1}= T_2^\vee \qquad\hbox{and}\qquad
g^\vee T_2^\vee (g^\vee)^{-1} = T_0^\vee.
$$
\end{prop}

To give a concrete example of Theorem \ref{Pexpansion}
let us compute the symmetric Macdonald polynomial
$P_\rho$ where $\rho = \alpha_1+\alpha_2$.
Since
$$
\mathbf{1_0} = X^{s_1s_2s_1} + t^{-1/2}X^{s_1s_2} +
t^{-1/2}X^{s_2s_1} + t^{-2/2}X^{s_1} + t^{-2/2}X^{s_2} + t^{-3/2},
$$
and $X^\rho m = s_0^\vee$ is the minimal length element
of the coset $X^\rho W_0$,
\begin{align*}
P_\rho &= \mathbf{1}_0E_\rho = \mathbf{1}_0\tau_0^\vee \mathbf{1}\\
    &= \left(X^{s_1s_2s_1} + t^{-1/2}X^{s_1s_2} + t^{-1/2}X^{s_2s_1} \right)
    \left(T_0^\vee + \frac{t^{-1/2}(1-t)}{1-Y^{-\alpha_0^\vee}}
    \right)\mathbf{1} \\
    &\qquad +
    \left(t^{-2/2}X^{s_1} + t^{-2/2}X^{s_2} + t^{-3/2} \right)
    \left((T_0^\vee)^{-1} +
    \frac{t^{-1/2}(1-t)Y^{-\alpha_0^\vee}}{1-Y^{-\alpha_0^\vee}}\right)
    \mathbf{1} \\
    &= \left(X^{s_1s_2s_1s_0} + t^{-1/2}X^{s_1s_2s_0} +
    t^{-1/2}X^{s_2s_1s_0} + t^{-2/2}X^{s_1s_0} + t^{-2/2}X^{s_2s_0} +
    t^{-3/2}X^{s_0}\right)\mathbf{1}\\
    &\qquad +\left(X^{s_1s_2s_1} + t^{-1/2}X^{s_1s_2} + t^{-1/2}X^{s_2s_1}\right)
    \frac{t^{-1/2}(1-t)}{1-Y^{-\alpha_0^\vee}}
    \mathbf{1}\\
    &\qquad +\left(t^{-2/2}X^{s_1} + t^{-2/2}X^{s_2} +
    t^{-3/2}\right)
    \frac{t^{-1/2}(1-t)Y^{-\alpha_0^\vee}}{1-Y^{-\alpha_0^\vee}}
    \mathbf{1}.
\end{align*}
Since
$Y^{-\alpha_0^\vee}\mathbf{1}
=Y^{\varphi^\vee-d}\mathbf{1} = qY^{\alpha^\vee_1 +\alpha^\vee_2}\mathbf{1}
= t^2q\mathbf{1}$,
\begin{align*}
P_\rho
    &= \left(
    \begin{array}{l}
    X^{w_0\rho} + t^{-1/2}X^{s_1s_2\rho}T_2^\vee +
    t^{-1/2}X^{s_2s_1\rho}T_1^\vee \\
    \qquad + t^{-2/2}X^{s_1\rho}T_2^\vee T_1^\vee
    + t^{-2/2}X^{s_2\rho}T_1^\vee T_2^\vee +
    t^{-3/2}X^{\rho}T_1^\vee T_2^\vee T_1^\vee
    \end{array} \right)\mathbf{1}\\
    &\qquad +\left(T_1^\vee T_2^\vee T_1^\vee + t^{-1/2}T_1^\vee T_2^\vee
    + t^{-1/2} T_2^\vee T_1^\vee\right)
    \frac{t^{-1/2}(1-t)}{1-t^2q}
    \mathbf{1}\\
    &\qquad +\left(t^{-2/2}T_1^\vee + t^{-2/2}T_2^\vee +
    t^{-3/2}\right)
    \frac{t^{-1/2}(1-t)t^2q}{1-t^2q}
    \mathbf{1}\\
&= \left(X^{w_0\rho} + X^{s_1s_2\rho} +
    X^{s_2s_1\rho} + X^{s_1\rho} + X^{s_2\rho} +
    X^{\rho}\right)\mathbf{1}\\
    &\qquad +\left(t^{3/2} + t^{1/2} + t^{1/2} \right)
    \frac{t^{-1/2}(1-t)}{1-t^2q}
    \mathbf{1} + \left(t^{-1/2} + t^{-1/2} + t^{-3/2}\right)
    \frac{t^{-1/2}(1-t)t^2q}{1-t^2q}
    \mathbf{1}\\
&= \left( X^{w_0\rho} + X^{s_1s_2\rho} +
    X^{s_2s_1\rho} + X^{s_1\rho} + X^{s_2\rho} +
    X^{\rho} + (t+2+2tq+q) \frac{1-t}{1-t^2q}
\right)\mathbf{1}.
\end{align*}
The set $\cP(\vec \rho)$ contains 12 alcove walks,
$$
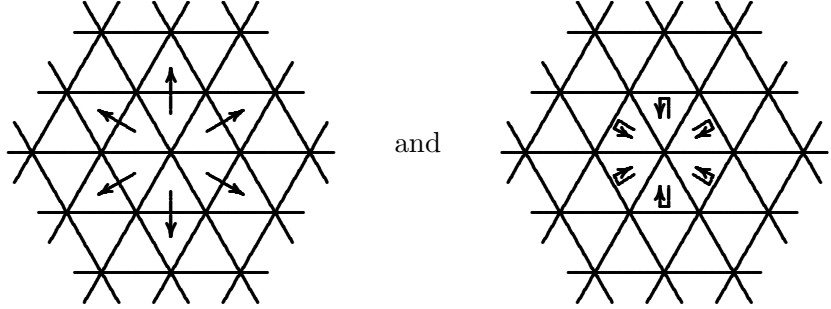
\beginpicture
\setcoordinatesystem units <0.4cm,0.4cm>         
\setplotarea x from -5.5 to 5.5, y from -5 to 5    
    \plot 1.732 -5  5.196 1 /
    \plot -0.577 -5  4.041 3 /
    \plot -2.887 -5  2.887 5 /
    \plot -4.041 -3  0.577 5 /
    \plot -5.196 -1  -1.732 5 /
    \plot -1.732 -5  -5.196 1 /
    \plot 0.577 -5  -4.041 3 /
    \plot 2.887 -5  -2.887 5 /
    \plot 4.041 -3  -0.577 5 /
    \plot 5.196 -1  1.732 5 /
    \plot -3.2 4  3.2 4 /
    \plot -4.4 2  4.4 2 /
    \plot -5.4 0  5.4 0 /
    \plot -3.2 -4  3.2 -4 /
    \plot -4.4 -2  4.4 -2 /
    \arrow <5pt> [.2,.67] from 0 1.3 to 0 2.8
    \arrow <5pt> [.2,.67] from 0 -1.3 to 0 -2.8
    \arrow <5pt> [.2,.67] from 1.2 0.7 to 2.4 1.4
    \arrow <5pt> [.2,.67] from 1.2 -0.7 to 2.4 -1.4
    \arrow <5pt> [.2,.67] from -1.2 0.7 to -2.4 1.4
    \arrow <5pt> [.2,.67] from -1.2 -0.7 to -2.4 -1.4
\endpicture\qquad\textrm{and}\qquad\beginpicture
\setcoordinatesystem units <0.4cm,0.4cm>         
\setplotarea x from -5.5 to 5.5, y from -5 to 5    
    \plot 1.732 -5  5.196 1 /
    \plot -0.577 -5  4.041 3 /
    \plot -2.887 -5  2.887 5 /
    \plot -4.041 -3  0.577 5 /
    \plot -5.196 -1  -1.732 5 /
    \plot -1.732 -5  -5.196 1 /
    \plot 0.577 -5  -4.041 3 /
    \plot 2.887 -5  -2.887 5 /
    \plot 4.041 -3  -0.577 5 /
    \plot 5.196 -1  1.732 5 /
    \plot -3.2 4  3.2 4 /
    \plot -4.4 2  4.4 2 /
    \plot -5.4 0  5.4 0 /
    \plot -3.2 -4  3.2 -4 /
    \plot -4.4 -2  4.4 -2 /
    \plot 0.15 -1.15 0.15 -1.8 /
    \plot 0.15 -1.8 -0.15 -1.8 /
    \arrow <5pt> [.2,.67] from -.15 -1.8 to -.15 -1.15
    \plot 0.15 1.15 0.15 1.8 /
    \plot 0.15 1.8 -0.15 1.8 /
    \arrow <5pt> [.2,.67] from -.15 1.8 to -.15 1.15
    \put{
\beginpicture
\setcoordinatesystem units <0.75cm,0.75cm>
\setplotarea x from -1 to 1, y from -1 to 1  
    \arrow <5pt> [.2,.67] from 0.3 -0.15 to 0 0
    \plot 0.3 -0.15 0.225 -0.3 /
    \plot 0.225 -0.3 -0.05 -0.13 /
\endpicture} at 0.9 -0.5
    \put{\beginpicture 
\setcoordinatesystem units <0.75cm,0.75cm>
\setplotarea x from -1 to 1, y from -1 to 1  
    \plot 0.05 -0.13 -0.225 -0.3 /
    \plot -0.225 -0.3 -0.3 -0.15 /
    \arrow <5pt> [.2,.67] from -0.3 -0.15 to 0 0
\endpicture} at -1.1 -0.5
    \put{
\beginpicture
\setcoordinatesystem units <0.75cm,0.75cm>
\setplotarea x from -1 to 1, y from -1 to 1  
    \arrow <5pt> [.2,.67] from 0.3 0.15 to 0 0
    \plot 0.3 0.15 0.225 0.3 /
    \plot 0.225 0.3 -0.05 0.13 /
\endpicture} at 0.9 0.5
    \put{\beginpicture 
\setcoordinatesystem units <0.75cm,0.75cm>
\setplotarea x from -1 to 1, y from -1 to 1  
    \plot 0.05 0.13 -0.225 0.3 /
    \plot -0.225 0.3 -0.3 0.15 /
    \arrow <5pt> [.2,.67] from -0.3 0.15 to 0 0
\endpicture} at -1.1 0.5
\endpicture$$
The Hall-Littlewood polynomial and the Weyl character are
$$P_\rho(0,t) = m_\rho+(2+t)(1-t)
\qquad\hbox{and}\qquad
s_\rho = P_\rho(0,0) = m_\rho+2,$$
where $m_\rho = X^{w_0\rho} + X^{s_1s_2\rho} +
    X^{s_2s_1\rho} + X^{s_1\rho} + X^{s_2\rho} +
    X^{\rho}$.

The expression $X^{s_1s_2\rho}s_2=s_1^\vee s_2^\vee s_0^\vee$
is a reduced word for the
minimal length element in the coset $X^{s_1s_2\rho}W_0$ and
Theorem \ref{mainthm} is illustrated by
\begin{align*}
\tau_1^\vee &\tau_2^\vee\tau_0^\vee
    = \left(T_1^\vee + \frac{t^{-1/2}(1-t)}{1-Y^{-\alpha_1^\vee}}
    \right) \tau_1^\vee \tau_0^\vee
= \left(T_1^\vee \tau_2^\vee +
    \tau_2^\vee\frac{t^{-1/2}(1-t)}{1-Y^{-s_2\alpha_1^\vee}}
    \right)\tau_0^\vee \\
&= \left(T_1^\vee T_2^\vee +
    T_1^\vee\frac{t^{-1/2}(1-t)}{1-Y^{-\alpha_2^\vee}} + 
    T_2^\vee\frac{t^{-1/2}(1-t)}{1-Y^{-s_2\alpha_1^\vee}}
    + \frac{t^{-1/2}(1-t)}{1-Y^{-\alpha_2^\vee}}\frac{t^{-1/2}(1-t)}{1-Y^{-s_2\alpha_1^\vee}}
    \right)\tau_0^\vee \\
&= T_1^\vee T_2^\vee \tau_0^\vee +
    T_1^\vee\tau_0^\vee\frac{t^{-1/2}(1-t)}{1-Y^{-s_0\alpha_2^\vee}} +
    T_2^\vee\tau_0^\vee\frac{t^{-1/2}(1-t)}{1-Y^{-s_0s_2\alpha_1^\vee}}
    + \tau_0^\vee\frac{t^{-1/2}(1-t)}{1-Y^{-s_0\alpha_2^\vee}}
    \frac{t^{-1/2}(1-t)}{1-Y^{-s_0s_2\alpha_1^\vee}}  \\
&= T_1^\vee T_2^\vee T_0^\vee + T_1^\vee T_2^\vee
    \frac{t^{-1/2}(1-t)}{1-Y^{-\alpha_0^\vee}} +
    T_1^\vee
    (T_0^\vee)^{-1}\frac{t^{-1/2}(1-t)}{1-Y^{-s_0\alpha_2^\vee}}  \\
    &\qquad\qquad\quad +
    T_1^\vee \frac{t^{-1/2}(1-t)Y^{-\alpha_0^\vee}}{1-Y^{-\alpha_0^\vee}}
    \frac{t^{-1/2}(1-t)}{1-Y^{-s_0\alpha_2^\vee}} +
    T_2^\vee (T_0^\vee)^{-1} \frac{t^{-1/2}(1-t)}{1-Y^{-s_0s_2\alpha_1^\vee}} \\
    &\qquad\qquad\quad
    + T_2^\vee \frac{t^{-1/2}(1-t)Y^{-\alpha_0^\vee}}{1-Y^{-\alpha_0^\vee}}
    \frac{t^{-1/2}(1-t)}{1-Y^{-s_0s_2\alpha_1^\vee}}
    + (T_0^\vee)^{-1}\frac{t^{-1/2}(1-t)}{1-Y^{-s_0\alpha_2^\vee}}
    \frac{t^{-1/2}(1-t)}{1-Y^{-s_0s_2\alpha_1^\vee}}  \\
    &\qquad\qquad\quad
    + \frac{t^{-1/2}(1-t)Y^{-\alpha_0^\vee}}{1-Y^{-\alpha_0^\vee}}
    \frac{t^{-1/2}(1-t)}{1-Y^{-s_0\alpha_2^\vee}}
    \frac{t^{-1/2}(1-t)}{1-Y^{-s_0s_2\alpha_1^\vee}},
\end{align*}
where the eight terms in this expansion correspond to the eight alcove walks in
$\cB(1,s_1^\vee s_2^\vee s_0^\vee)=\cB(\overrightarrow{s_1s_2\rho})$ pictured below.
Applying the expansion of $\tau_1^\vee \tau_2^\vee \tau_0^\vee$
to $\mathbf{1}$ and using
\begin{equation}
\begin{array}{c}
Y^{-\alpha_0^\vee}\mathbf{1}
= Y^{\varphi^\vee-d}\mathbf{1} = t^2q\mathbf{1},
\qquad\quad
Y^{-s_0\alpha_2^\vee} \mathbf{1}
= Y^{\alpha_1^\vee -d}\mathbf{1} = tq\mathbf{1},
\\ \\
\hbox{and}\qquad
Y^{-s_0s_2\alpha_1^\vee} \mathbf{1}
=Y^{\varphi^\vee-2d}\mathbf{1} = t^2q^2 \mathbf{1},
\end{array}
\end{equation}
computes
\begin{align*}
E_{s_1s_2\rho} &= \left(X^{s_1s_2\rho}t^{1/2}
+ t \frac{t^{-1/2}(1-t)}{1-t^2q}
+  X^{s_1\rho} t\frac{t^{-1/2}(1-t)}{1-tq}
+ t^{1/2} \frac{t^{-1/2}(1-t)t^2q}{1-t^2q} \frac{t^{-1/2}(1-t)}{1-tq} \right. \\
&\qquad\qquad\quad
+ X^{s_2\rho} t\frac{t^{-1/2}(1-t)}{1-t^2q^2}
+ t^{1/2} \frac{t^{-1/2}(1-t)t^2q}{1-t^2q} \frac{t^{-1/2}(1-t)}{1-t^2q^2} \\
&\qquad\qquad\quad \left.
+ X^\rho t^{3/2}\frac{t^{-1/2}(1-t)}{1-tq} \frac{t^{-1/2}(1-t)}{1-t^2q^2}
+ \frac{t^{-1/2}(1-t)t^2q}{1-t^2q}\frac{t^{-1/2}(1-t)}{1-tq}
  \frac{t^{-1/2}(1-t)}{1-t^2q^2} \right)
\mathbf{1}. \\
&= t^{1/2}\left(X^{s_1s_2\rho}
+  \frac{(1-t)}{1-t^2q}
+  X^{s_1\rho} \frac{(1-t)}{1-tq}
+ t \frac{(1-t)q}{1-t^2q} \frac{(1-t)}{1-tq}
+ X^{s_2\rho} \frac{(1-t)}{1-t^2q^2} \right. \\
&\qquad\qquad\quad\left.
+ t \frac{(1-t)q}{1-t^2q} \frac{(1-t)}{1-t^2q^2} 
+ X^\rho \frac{(1-t)}{1-tq} \frac{(1-t)}{1-t^2q^2}
+ \frac{(1-t)q}{1-t^2q}\frac{(1-t)}{1-tq}
  \frac{(1-t)}{1-t^2q^2} \right)
\mathbf{1}.
\end{align*}

$$\begin{matrix}
\beginpicture
\setcoordinatesystem units <0.4cm,0.4cm>         
\setplotarea x from -5.5 to 5.5, y from -5 to 5    
    \plot 1.732 -5  5.196 1 /
    \plot -0.577 -5  4.041 3 /
    \plot -2.887 -5  2.887 5 /
    \plot -4.041 -3  0.577 5 /
    \plot -5.196 -1  -1.732 5 /
    \plot -1.732 -5  -5.196 1 /
    \plot 0.577 -5  -4.041 3 /
    \plot 2.887 -5  -2.887 5 /
    \plot 4.041 -3  -0.577 5 /
    \plot 5.196 -1  1.732 5 /
    \plot -3.2 4  3.2 4 /
    \plot -4.4 2  4.4 2 /
    \plot -5.4 0  5.4 0 /
    \plot -3.2 -4  3.2 -4 /
    \plot -4.4 -2  4.4 -2 /
    \arrow <5pt> [.2,.67] from 1.2 0.7 to 1.2 -0.7 
    \arrow <5pt> [.2,.67] from 1.2 -0.7 to 2.4 -1.4 
    \arrow <5pt> [.2,.67] from 0 1.3 to 1.2 0.7 
\put{$\bullet$} at 0 0
\endpicture
&
\beginpicture
\setcoordinatesystem units <0.4cm,0.4cm>         
\setplotarea x from -5.5 to 5.5, y from -5 to 5    
    \plot 1.732 -5  5.196 1 /
    \plot -0.577 -5  4.041 3 /
    \plot -2.887 -5  2.887 5 /
    \plot -4.041 -3  0.577 5 /
    \plot -5.196 -1  -1.732 5 /
    \plot -1.732 -5  -5.196 1 /
    \plot 0.577 -5  -4.041 3 /
    \plot 2.887 -5  -2.887 5 /
    \plot 4.041 -3  -0.577 5 /
    \plot 5.196 -1  1.732 5 /
    \plot -3.2 4  3.2 4 /
    \plot -4.4 2  4.4 2 /
    \plot -5.4 0  5.4 0 /
    \plot -3.2 -4  3.2 -4 /
    \plot -4.4 -2  4.4 -2 /
    \arrow <5pt> [.2,.67] from 0 1.3 to 1.2 0.7
    \arrow <5pt> [.2,.67] from 1.2 0.7 to 1.2 -0.7
    \put{\beginpicture
    \setcoordinatesystem units <0.75cm,0.75cm>
    \setplotarea x from -1 to 1, y from -1 to 1  
    \arrow <5pt> [.2,.67] from 0.3 -0.15 to 0 0
    \plot 0.3 -0.15 0.225 -0.3 /
    \plot 0.225 -0.3 -0.05 -0.13 / 
    \endpicture} at 1.3 -0.5
    \put{$\bullet$} at 0 0
\endpicture
&
\beginpicture
\setcoordinatesystem units <0.4cm,0.4cm>         
\setplotarea x from -5.5 to 5.5, y from -5 to 5    
    \plot 1.732 -5  5.196 1 /
    \plot -0.577 -5  4.041 3 /
    \plot -2.887 -5  2.887 5 /
    \plot -4.041 -3  0.577 5 /
    \plot -5.196 -1  -1.732 5 /
    \plot -1.732 -5  -5.196 1 /
    \plot 0.577 -5  -4.041 3 /
    \plot 2.887 -5  -2.887 5 /
    \plot 4.041 -3  -0.577 5 /
    \plot 5.196 -1  1.732 5 /
    \plot -3.2 4  3.2 4 /
    \plot -4.4 2  4.4 2 /
    \plot -5.4 0  5.4 0 /
    \plot -3.2 -4  3.2 -4 /
    \plot -4.4 -2  4.4 -2 /
    \arrow <5pt> [.2,.67] from 0 1.3 to 1.2 0.7
    \arrow <5pt> [.2,.67] from 1.5 0.75 to 2.4 1.4
    \plot 1.2 0.65 1.2 0 /
    \plot 1.2 0 1.5 0   /
    \arrow <5pt> [.2,.67] from 1.5 0 to 1.5 0.75
    \put{$\bullet$} at 0 0
\endpicture
\\ \\
X^{s_2s_1\rho}t^{1/2}
&\displaystyle{ t^{1/2}\frac{(1-t)}{1-t^2q} }
&\displaystyle{ X^{s_1\rho} t^{1/2} \frac{(1-t)}{1-tq}}
\end{matrix}
$$
$$
\begin{matrix}
\beginpicture
\setcoordinatesystem units <0.4cm,0.4cm>         
\setplotarea x from -5.5 to 5.5, y from -5 to 5    
    \plot 1.732 -5  5.196 1 /
    \plot -0.577 -5  4.041 3 /
    \plot -2.887 -5  2.887 5 /
    \plot -4.041 -3  0.577 5 /
    \plot -5.196 -1  -1.732 5 /
    \plot -1.732 -5  -5.196 1 /
    \plot 0.577 -5  -4.041 3 /
    \plot 2.887 -5  -2.887 5 /
    \plot 4.041 -3  -0.577 5 /
    \plot 5.196 -1  1.732 5 /
    \plot -3.2 4  3.2 4 /
    \plot -4.4 2  4.4 2 /
    \plot -5.4 0  5.4 0 /
    \plot -3.2 -4  3.2 -4 /
    \plot -4.4 -2  4.4 -2 /
    \arrow <5pt> [.2,.67] from 0 1.3 to 1.2 0.7
    \plot 1.2 0.65 1.2 0 /
    \plot 1.2 0 1.5 0   /
    \arrow <5pt> [.2,.67] from 1.5 0 to 1.5 0.75
    \put{
    \beginpicture
    \setcoordinatesystem units <0.75cm,0.75cm>
    \setplotarea x from -1 to 1, y from -1 to 1  
    \arrow <5pt> [.2,.67] from 0.3 0.15 to 0 0
    \plot 0.3 0.15 0.225 0.3 /
    \plot 0.225 0.3 -0.05 0.13 /
    \endpicture} at 1.4 0.55
    \put{$\bullet$} at 0 0
\endpicture
&
\beginpicture
\setcoordinatesystem units <0.4cm,0.4cm>         
\setplotarea x from -5.5 to 5.5, y from -5 to 5    
    \plot 1.732 -5  5.196 1 /
    \plot -0.577 -5  4.041 3 /
    \plot -2.887 -5  2.887 5 /
    \plot -4.041 -3  0.577 5 /
    \plot -5.196 -1  -1.732 5 /
    \plot -1.732 -5  -5.196 1 /
    \plot 0.577 -5  -4.041 3 /
    \plot 2.887 -5  -2.887 5 /
    \plot 4.041 -3  -0.577 5 /
    \plot 5.196 -1  1.732 5 /
    \plot -3.2 4  3.2 4 /
    \plot -4.4 2  4.4 2 /
    \plot -5.4 0  5.4 0 /
    \plot -3.2 -4  3.2 -4 /
    \plot -4.4 -2  4.4 -2 /
    \put{
    \beginpicture
    \setcoordinatesystem units <0.75cm,0.75cm>
    \setplotarea x from -1 to 1, y from -1 to 1  
    \arrow <5pt> [.2,.67] from 0.3 -0.15 to 0 0
    \plot 0.3 -0.15 0.225 -0.3 /
    \plot 0.225 -0.3 -0.05 -0.13 /
    \endpicture} at -0.15 1.3
    \arrow <5pt> [.2,.67] from 0 1.3 to -1.2 0.7
    \arrow <5pt> [.2,.67] from -1.2 0.7 to -2.4 1.4
    \put{$\bullet$} at 0 0
\endpicture
&
\beginpicture
\setcoordinatesystem units <0.4cm,0.4cm>         
\setplotarea x from -5.5 to 5.5, y from -5 to 5    
    \plot 1.732 -5  5.196 1 /
    \plot -0.577 -5  4.041 3 /
    \plot -2.887 -5  2.887 5 /
    \plot -4.041 -3  0.577 5 /
    \plot -5.196 -1  -1.732 5 /
    \plot -1.732 -5  -5.196 1 /
    \plot 0.577 -5  -4.041 3 /
    \plot 2.887 -5  -2.887 5 /
    \plot 4.041 -3  -0.577 5 /
    \plot 5.196 -1  1.732 5 /
    \plot -3.2 4  3.2 4 /
    \plot -4.4 2  4.4 2 /
    \plot -5.4 0  5.4 0 /
    \plot -3.2 -4  3.2 -4 /
    \plot -4.4 -2  4.4 -2 /
    \put{
    \beginpicture
    \setcoordinatesystem units <0.75cm,0.75cm>
    \setplotarea x from -1 to 1, y from -1 to 1  
    \arrow <5pt> [.2,.67] from 0.3 -0.15 to 0 0
    \plot 0.3 -0.15 0.225 -0.3 /
    \plot 0.225 -0.3 -0.05 -0.13 /
    \endpicture} at -0.15 1.3
    \arrow <5pt> [.2,.67] from 0 1.3 to -1.2 0.7
    \put{\beginpicture 
    \setcoordinatesystem units <0.75cm,0.75cm>
    \setplotarea x from -1 to 1, y from -1 to 1  
    \plot 0.05 0.13 -0.225 0.3 /
    \plot -0.225 0.3 -0.3 0.15 /
    \arrow <5pt> [.2,.67] from -0.3 0.15 to 0 0
    \endpicture} at -1.35 0.4
    \put{$\bullet$} at 0 0
\endpicture
\\ \\
\displaystyle{ t^{3/2}\frac{(1-t)}{1-tq}\frac{(1-t)q}{1-t^2q} }
&\displaystyle{ X^{s_2\rho} t^{1/2} \frac{(1-t)}{1-t^2q^2}}
&\displaystyle{t^{3/2}\frac{(1-t)q}{1-t^2q}\frac{(1-t)}{1-t^2q^2} }
\end{matrix}
$$

$$\begin{matrix}
\beginpicture
\setcoordinatesystem units <0.4cm,0.4cm>         
\setplotarea x from -5.5 to 5.5, y from -5 to 5    
    \plot 1.732 -5  5.196 1 /
    \plot -0.577 -5  4.041 3 /
    \plot -2.887 -5  2.887 5 /
    \plot -4.041 -3  0.577 5 /
    \plot -5.196 -1  -1.732 5 /
    \plot -1.732 -5  -5.196 1 /
    \plot 0.577 -5  -4.041 3 /
    \plot 2.887 -5  -2.887 5 /
    \plot 4.041 -3  -0.577 5 /
    \plot 5.196 -1  1.732 5 /
    \plot -3.2 4  3.2 4 /
    \plot -4.4 2  4.4 2 /
    \plot -5.4 0  5.4 0 /
    \plot -3.2 -4  3.2 -4 /
    \plot -4.4 -2  4.4 -2 /
    \put{
    \beginpicture
    \setcoordinatesystem units <0.75cm,0.75cm>
    \setplotarea x from -1 to 1, y from -1 to 1  
    \arrow <5pt> [.2,.67] from 0.3 -0.15 to 0 0
    \plot 0.3 -0.15 0.225 -0.3 /
    \plot 0.225 -0.3 -0.05 -0.13 /
    \endpicture} at -0.15 1.3
    \put{\beginpicture 
    \setcoordinatesystem units <0.75cm,0.75cm>
    \setplotarea x from -1 to 1, y from -1 to 1  
    \plot 0.05 -0.13 -0.225 -0.3 /
    \plot -0.225 -0.3 -0.3 -0.15 /
    \arrow <5pt> [.2,.67] from -0.3 -0.15 to 0 0
    \endpicture} at -0.1 1.6
    \arrow <5pt> [.2,.67] from 0 1.6 to 0 2.8
    \put{$\bullet$} at 0 0
\endpicture
\qquad\qquad
&\qquad\qquad
\beginpicture
\setcoordinatesystem units <0.4cm,0.4cm>         
\setplotarea x from -5.5 to 5.5, y from -5 to 5    
    \plot 1.732 -5  5.196 1 /
    \plot -0.577 -5  4.041 3 /
    \plot -2.887 -5  2.887 5 /
    \plot -4.041 -3  0.577 5 /
    \plot -5.196 -1  -1.732 5 /
    \plot -1.732 -5  -5.196 1 /
    \plot 0.577 -5  -4.041 3 /
    \plot 2.887 -5  -2.887 5 /
    \plot 4.041 -3  -0.577 5 /
    \plot 5.196 -1  1.732 5 /
    \plot -3.2 4  3.2 4 /
    \plot -4.4 2  4.4 2 /
    \plot -5.4 0  5.4 0 /
    \plot -3.2 -4  3.2 -4 /
    \plot -4.4 -2  4.4 -2 /
    \put{
    \beginpicture
    \setcoordinatesystem units <0.75cm,0.75cm>
    \setplotarea x from -1 to 1, y from -1 to 1  
    \arrow <5pt> [.2,.67] from 0.3 -0.15 to 0 0
    \plot 0.3 -0.15 0.225 -0.3 /
    \plot 0.225 -0.3 -0.05 -0.13 /
    \endpicture} at -0.15 1.3
    \put{\beginpicture 
    \setcoordinatesystem units <0.75cm,0.75cm>
    \setplotarea x from -1 to 1, y from -1 to 1  
    \plot 0.05 -0.13 -0.225 -0.3 /
    \plot -0.225 -0.3 -0.3 -0.15 /
    \arrow <5pt> [.2,.67] from -0.3 -0.15 to 0 0
    \endpicture} at -0.1 1.6
    \plot -0.15 1.6 -0.15 2.2 /
    \plot 0.15 2.2 -0.15 2.2 /
    \arrow <5pt> [.2,.67] from .15 2.2 to .15 1.4
    \put{$\bullet$} at 0 0
\endpicture
\\ \\
\displaystyle{ X^\rho t^{1/2} \frac{(1-t)}{1-tq} \frac{(1-t)}{1-t^2q^2} }\qquad\quad
&\quad\qquad
\displaystyle{ t^{3/2}\frac{(1-t)}{1-tq}\frac{(1-t)q}{1-t^2q} \frac{(1-t)}{1-t^2q^2} }
\end{matrix}
$$

\newpage

\begin{section}{Appendix: The bijection between $W$ and alcoves in
type $SL_3$}

The following pictures illustrate the bijection of \eqref{Wtoalcoves} for type $SL_3$.  In this
case, $\Omega^\vee = \{1,g^\vee,(g^\vee)^2\}\cong \ZZ/3\ZZ$, and
$\Omega^\vee\times \fh_\RR^*$ has 3 sheets.
The alcoves are the triangles and the (centres of) hexagons are the elements of $\fh_\ZZ^*$.

$$
\beginpicture
\setcoordinatesystem units <1.4cm,1.4cm>         
\setplotarea x from -5 to 5, y from -5 to 5.2  
    \plot 2.8 -3.81 3.7 -2.2516 /
    \plot 1.8 -3.81 3.7 -0.52 /
    \plot 0.8 -3.81 3.7 1.212 /
    \plot -0.2 -3.81 3.7 3.044 /
    \plot -1.2 -3.81 3.2 3.81 /
    \plot -2.2 -3.81 2.2 3.81 /
    \plot -3.2 -3.81 1.2 3.81 /
    \plot -3.7 -3.044 0.2 3.81 /
    \plot -3.7 -1.212 -0.8 3.81 /
    \plot -3.7 0.520 -1.8 3.81 /
    \plot -3.7 2.2516 -2.8 3.81 /
    \plot -2.8 -3.81 -3.7 -2.2516 /
    \plot -1.8 -3.81 -3.7 -0.52 /
    \plot -0.8 -3.81 -3.7 1.212 /
    \plot 0.2 -3.81 -3.7 3.044 /
    \plot 1.2 -3.81 -3.2 3.81 /
    \plot 2.2 -3.81 -2.2 3.81 /
    \plot 3.2 -3.81 -1.2 3.81 /
    \plot 3.7 -3.044 -0.2 3.81 /
    \plot 3.7 -1.212 0.8 3.81 /
    \plot 3.7 0.520 1.8 3.81 /
    \plot 3.7 2.2516 2.8 3.81 /
    \plot -3.7 -3.464 3.7 -3.464 /
    \plot -3.7 -2.598 3.7 -2.598 /
    \plot -3.7 -1.732 3.7 -1.732 /
    \plot -3.7 -0.866 3.7 -0.866 /
    \plot -3.7 0 3.7 0 /
    \plot -3.7 3.464 3.7 3.464 /
    \plot -3.7 2.598 3.7 2.598 /
    \plot -3.7 1.732 3.7 1.732 /
    \plot -3.7 0.866 3.7 0.866 /
    \put{$\bullet$} at 0 3.464
    \put{$\bullet$} at 0 -3.464
    \put{$\bullet$} at 3 3.464
    \put{$\bullet$} at 3 -3.464
    \put{$\bullet$} at -3 3.464
    \put{$\bullet$} at -3 -3.464
    \small{\put{$1$} at 0 0.5
    \put{$s^\vee_0$} at 0 1.1
    \put{$s^\vee_1$} at 0.5 0.3
    \put{$s^\vee_2$} at -0.5 0.3
    \put{$s^\vee_0s^\vee_1$} at 0.5 1.55
    \put{$s^\vee_0s^\vee_2$} at -0.5 1.55
    \put{$s^\vee_1s^\vee_2$} at 0.5 -0.2
    \put{$s^\vee_2s^\vee_1$} at -0.5 -0.2
    \put{$s^\vee_2s^\vee_0$} at -1 0.7
    \put{$s^\vee_1s^\vee_0$} at 1 0.7
    \put{$w_0$} at 0 -0.65
    \put{$X^\rho$} at 0 2.4
    \put{$X^{s_1\rho}$} at 1.5 1.5
    \put{$X^{s_2\rho}$} at -1.5 1.5
    \put{$X^{s_1s_2\rho}$} at 1.5 -0.2
    \put{$X^{s_1\rho}$} at -1.5 -0.2
    \put{$X^{w_0\rho}$} at 0 -1.1 }%
    \put{\beginpicture
    \setplotsymbol({\tiny{$\bullet$}})
    \plot 1 0 0.5 0.866 /
    \plot 0.5 0.866 -0.5 0.866 /
    \plot -0.5 0.866 -1 0 /
    \plot -1 0 -0.5 -0.866 /
    \plot -0.5 -0.866 0.5 -0.866 /
    \plot 0.5 -0.866 1 0 /
    \put{$\bullet$} at 0 0
    \endpicture} at 0 0
    \put{\beginpicture
    \setplotsymbol({\tiny{$\bullet$}})
    \plot 1 0 0.5 0.866 /
    \plot 0.5 0.866 -0.5 0.866 /
    \plot -0.5 0.866 -1 0 /
    \plot -1 0 -0.5 -0.866 /
    \plot -0.5 -0.866 0.5 -0.866 /
    \plot 0.5 -0.866 1 0 /
    \put{$\bullet$} at 0 0
    \endpicture} at 0 1.732
    \put{\beginpicture
    \setplotsymbol({\tiny{$\bullet$}})
    \plot 1 0 0.5 0.866 /
    \plot 0.5 0.866 -0.5 0.866 /
    \plot -0.5 0.866 -1 0 /
    \plot -1 0 -0.5 -0.866 /
    \plot -0.5 -0.866 0.5 -0.866 /
    \plot 0.5 -0.866 1 0 /
    \put{$\bullet$} at 0 0
    \endpicture} at 0 -1.732
    \put{\beginpicture
    \setplotsymbol({\tiny{$\bullet$}})
    \plot 1 0 0.5 0.866 /
    \plot 0.5 0.866 -0.5 0.866 /
    \plot -0.5 0.866 -1 0 /
    \plot -1 0 -0.5 -0.866 /
    \plot -0.5 -0.866 0.5 -0.866 /
    \plot 0.5 -0.866 1 0 /
    \put{$\bullet$} at 0 0
    \endpicture} at -1.5 0.866
    \put{\beginpicture
    \setplotsymbol({\tiny{$\bullet$}})
    \plot 1 0 0.5 0.866 /
    \plot 0.5 0.866 -0.5 0.866 /
    \plot -0.5 0.866 -1 0 /
    \plot -1 0 -0.5 -0.866 /
    \plot -0.5 -0.866 0.5 -0.866 /
    \plot 0.5 -0.866 1 0 /
    \put{$\bullet$} at 0 0
    \endpicture} at -1.5 -0.866
    \put{\beginpicture
    \setplotsymbol({\tiny{$\bullet$}})
    \plot 1 0 0.5 0.866 /
    \plot 0.5 0.866 -0.5 0.866 /
    \plot -0.5 0.866 -1 0 /
    \plot -1 0 -0.5 -0.866 /
    \plot -0.5 -0.866 0.5 -0.866 /
    \plot 0.5 -0.866 1 0 /
    \put{$\bullet$} at 0 0
    \endpicture} at -1.5 2.598
    \put{\beginpicture
    \setplotsymbol({\tiny{$\bullet$}})
    \plot 1 0 0.5 0.866 /
    \plot 0.5 0.866 -0.5 0.866 /
    \plot -0.5 0.866 -1 0 /
    \plot -1 0 -0.5 -0.866 /
    \plot -0.5 -0.866 0.5 -0.866 /
    \plot 0.5 -0.866 1 0 /
    \put{$\bullet$} at 0 0
    \endpicture} at -1.5 -2.598
    \put{\beginpicture
    \setplotsymbol({\tiny{$\bullet$}})
    \plot 1 0 0.5 0.866 /
    \plot 0.5 0.866 -0.5 0.866 /
    \plot -0.5 0.866 -1 0 /
    \plot -1 0 -0.5 -0.866 /
    \plot -0.5 -0.866 0.5 -0.866 /
    \plot 0.5 -0.866 1 0 /
    \put{$\bullet$} at 0 0
    \endpicture} at 1.5 0.866
    \put{\beginpicture
    \setplotsymbol({\tiny{$\bullet$}})
    \plot 1 0 0.5 0.866 /
    \plot 0.5 0.866 -0.5 0.866 /
    \plot -0.5 0.866 -1 0 /
    \plot -1 0 -0.5 -0.866 /
    \plot -0.5 -0.866 0.5 -0.866 /
    \plot 0.5 -0.866 1 0 /
    \put{$\bullet$} at 0 0
    \endpicture} at 1.5 -0.866
    \put{\beginpicture
    \setplotsymbol({\tiny{$\bullet$}})
    \plot 1 0 0.5 0.866 /
    \plot 0.5 0.866 -0.5 0.866 /
    \plot -0.5 0.866 -1 0 /
    \plot -1 0 -0.5 -0.866 /
    \plot -0.5 -0.866 0.5 -0.866 /
    \plot 0.5 -0.866 1 0 /
    \put{$\bullet$} at 0 0
    \endpicture} at 1.5 2.598
    \put{\beginpicture
    \setplotsymbol({\tiny{$\bullet$}})
    \plot 1 0 0.5 0.866 /
    \plot 0.5 0.866 -0.5 0.866 /
    \plot -0.5 0.866 -1 0 /
    \plot -1 0 -0.5 -0.866 /
    \plot -0.5 -0.866 0.5 -0.866 /
    \plot 0.5 -0.866 1 0 /
    \put{$\bullet$} at 0 0
    \endpicture} at 1.5 -2.598
    \put{\beginpicture
    \setplotsymbol({\tiny{$\bullet$}})
    \plot 0.5 -0.866 -0.5 -0.866 /
    \plot -0.5 -0.866 -1 0 /
    \plot -1 0 -0.5 0.866 /
    \plot -0.5 0.866 0.5 0.866 /
    \put{$\bullet$} at 0 0
    \endpicture} at 3 -1.732
    \put{\beginpicture
    \setplotsymbol({\tiny{$\bullet$}})
    \plot 0.5 -0.866 -0.5 -0.866 /
    \plot -0.5 -0.866 -1 0 /
    \plot -1 0 -0.5 0.866 /
    \plot -0.5 0.866 0.5 0.866 /
    \put{$\bullet$} at 0 0
    \endpicture} at 3 0
    \put{\beginpicture
    \setplotsymbol({\tiny{$\bullet$}})
    \plot 0.5 -0.866 -0.5 -0.866 /
    \plot -0.5 0.866 0.5 0.866 /
    \put{$\bullet$} at 0 0
    \endpicture} at 3 1.732
    \put{\beginpicture
    \setplotsymbol({\tiny{$\bullet$}})
    \plot 0.5 -0.866 -0.5 -0.866 /
    \plot -0.5 0.866 0.5 0.866 /
    \put{$\bullet$} at 0 0
    \endpicture} at -3 1.732
    \put{\beginpicture
    \setplotsymbol({\tiny{$\bullet$}})
    \plot 0.5 -0.866 -0.5 -0.866 /
    \plot -0.5 0.866 0.5 0.866 /
    \put{$\bullet$} at 0 0
    \endpicture} at -3 0
    \put{\beginpicture
    \setplotsymbol({\tiny{$\bullet$}})
    \plot 0.5 -0.866 -0.5 -0.866 /
    \plot -0.5 0.866 0.5 0.866 /
    \put{$\bullet$} at 0 0
    \endpicture} at -3 -1.732
    \plot 2.5 -4.33 3.7 -2.2516 /
    \plot 1.5 -4.33 3.7 -0.52 /
    \plot 0.5 -4.33 3.7 1.212 /
    \plot -0.5 -4.33 3.7 3.044 /
    \plot -1.5 -4.33 3.2 3.81 /
    \plot -2.5 -4.33 2.2 3.81 /
    \plot -3.5 -4.33 1.2 3.81 /
    \plot -3.7 -3.044 0.2 3.81 /
    \plot -3.7 -1.212 -0.8 3.81 /
    \plot -3.7 0.520 -1.8 3.81 /
    \plot -3.7 2.2516 -2.8 3.81 /
    \plot 2.5 4.33 3.7 2.2516 /
    \plot 1.5 4.33 3.7 0.52 /
    \plot 0.5 4.33 3.7 -1.212 /
    \plot -0.5 4.33 3.7 -3.044 /
    \plot -1.5 4.33 3.2 -3.81 /
    \plot -2.5 4.33 2.2 -3.81 /
    \plot -3.5 4.33 1.2 -3.81 /
    \plot -3.7 3.044 0.2 -3.81 /
    \plot -3.7 1.212 -0.8 -3.81 /
    \plot -3.7 -0.520 -1.8 -3.81 /
    \plot -3.7 -2.2516 -2.8 -3.81 /
    \plot -3.7 -3.464 4.5 -3.464 /
    \plot -3.7 -2.598 4.5 -2.598 /
    \plot -3.7 -1.732 4.5 -1.732 /
    \plot -3.7 -0.866 4.5 -0.866 /
    \plot -3.7 0 4.5 0 /
    \plot -3.7 3.464 4.5 3.464 /
    \plot -3.7 2.598 4.5 2.598 /
    \plot -3.7 1.732 4.5 1.732 /
    \plot -3.7 0.866 4.5 0.866 /
    {\small\put{$\fh^{-\alpha^\vee_1+d}$} at -3.5 -4.7
    \put{$\fh^{\alpha_1^\vee}$} at -2.5 -4.7
    \put{$\fh^{\alpha_1^\vee+2d}$} at -0.5 -4.7
    \put{$\fh^{\alpha_1^\vee+4d}$} at 1.5 -4.7 } %
    {\small\put{$\fh^{\alpha_2^\vee+d}$} at -3.5 4.7
    \put{$\fh^{\alpha_2^\vee}$} at -2.5 4.7
    \put{$\fh^{-\alpha_2^\vee+2d}$} at -0.5 4.7
    \put{$\fh^{-\alpha_2^\vee+4d}$} at 1.5 4.7 } %
    {\small\put{$\fh^{-\varphi+4d}$}[l] at 4.5 3.464
    \put{$\fh^{-\varphi^\vee+3d}$}[l] at 4.5 2.598
    \put{$\fh^{-\varphi^\vee+2d}$}[l] at 4.5 1.732
    \put{$\fh^{\alpha_0^\vee}$}[l] at 4.5 0.866
    \put{$\fh^{\varphi^\vee}$}[l] at 4.5 0
    \put{$\fh^{\varphi^\vee+d}$}[l] at 4.5 -0.866
    \put{$\fh^{\varphi^\vee+2d}$}[l] at 4.5 -1.732
    \put{$\fh^{\varphi^\vee+3d}$}[l] at 4.5 -2.598
    \put{$\fh^{\varphi^\vee+4d}$}[l] at 4.5 -3.464 }%
    \put{$+$} at 4.25 0.2
    \put{$+$} at 4.25 1.066
    \put{$+$} at 4.25 1.932
    \put{$+$} at 4.25 2.798
    \put{$+$} at 4.25 3.664
    \put{$-$} at 4.25 -1.066
    \put{$-$} at 4.25 -1.932
    \put{$-$} at 4.25 -2.798
    \put{$-$} at 4.25 -3.664
    \put{$+$} at 4.25 -3.264
    \put{$+$} at 4.25 -2.398
    \put{$+$} at 4.25 -1.532
    \put{$+$} at 4.25 -0.666
    \put{$-$} at 4.25 3.264
    \put{$-$} at 4.25 2.398
    \put{$-$} at 4.25 1.532
    \put{$-$} at 4.25 0.666
    \put{$-$} at 4.25 -0.2
    \put{$+$} at -0.9 -4.2
    \put{$-$} at -0.4 -4.2
    \put{$+$} at -1.9 -4.2
    \put{$-$} at -1.4 -4.2
    \put{$+$} at -2.9 -4.2
    \put{$-$} at -2.4 -4.2
    \put{$+$} at -3.9 -4.2
    \put{$-$} at -3.4 -4.2
    \put{$+$} at 0.1 -4.2
    \put{$-$} at 0.6 -4.2
    \put{$+$} at 1.1 -4.2
    \put{$-$} at 1.6 -4.2
    \put{$+$} at 2.1 -4.2
    \put{$-$} at 2.6 -4.2
    \put{$-$} at -0.9 4.2
    \put{$+$} at -0.4 4.2
    \put{$-$} at -1.9 4.2
    \put{$+$} at -1.4 4.2
    \put{$-$} at -2.9 4.2
    \put{$+$} at -2.4 4.2
    \put{$-$} at -3.9 4.2
    \put{$+$} at -3.4 4.2
    \put{$-$} at 0.1 4.2
    \put{$+$} at 0.6 4.2
    \put{$-$} at 1.1 4.2
    \put{$+$} at 1.6 4.2
    \put{$-$} at 2.1 4.2
    \put{$+$} at 2.6 4.2
    \put{Sheet $1$} at 0 -5.3
\endpicture
$$

$$\beginpicture
\setcoordinatesystem units <1.4cm,1.4cm>         
\setplotarea x from -5 to 5, y from -5.2 to 5.2  
    \plot 2.8 -3.81 3.7 -2.2516 /
    \plot 1.8 -3.81 3.7 -0.52 /
    \plot 0.8 -3.81 3.7 1.212 /
    \plot -0.2 -3.81 3.7 3.044 /
    \plot -1.2 -3.81 3.2 3.81 /
    \plot -2.2 -3.81 2.2 3.81 /
    \plot -3.2 -3.81 1.2 3.81 /
    \plot -3.7 -3.044 0.2 3.81 /
    \plot -3.7 -1.212 -0.8 3.81 /
    \plot -3.7 0.520 -1.8 3.81 /
    \plot -3.7 2.2516 -2.8 3.81 /
    \plot -2.8 -3.81 -3.7 -2.2516 /
    \plot -1.8 -3.81 -3.7 -0.52 /
    \plot -0.8 -3.81 -3.7 1.212 /
    \plot 0.2 -3.81 -3.7 3.044 /
    \plot 1.2 -3.81 -3.2 3.81 /
    \plot 2.2 -3.81 -2.2 3.81 /
    \plot 3.2 -3.81 -1.2 3.81 /
    \plot 3.7 -3.044 -0.2 3.81 /
    \plot 3.7 -1.212 0.8 3.81 /
    \plot 3.7 0.520 1.8 3.81 /
    \plot 3.7 2.2516 2.8 3.81 /
    \plot -3.7 -3.464 3.7 -3.464 /
    \plot -3.7 -2.598 3.7 -2.598 /
    \plot -3.7 -1.732 3.7 -1.732 /
    \plot -3.7 -0.866 3.7 -0.866 /
    \plot -3.7 0 3.7 0 /
    \plot -3.7 3.464 3.7 3.464 /
    \plot -3.7 2.598 3.7 2.598 /
    \plot -3.7 1.732 3.7 1.732 /
    \plot -3.7 0.866 3.7 0.866 /
    \put{$\bullet$} at 0 0
    \put{$\bullet$} at -2 3.464
    \put{$\bullet$} at -2 -3.464
    \put{$\bullet$} at 2.5 3.464
    \put{$\bullet$} at 2.5 -3.464
    \put{$\bullet$} at -3.5 0.866
    \put{$\bullet$} at -3.5 0.866
    \put{$\bullet$} at -3.5 2.598
    \put{$\bullet$} at -3.5 -2.598
    \small{\put{$g^\vee$} at 0 0.55
    \put{$g^\vee s^\vee_2$} at 0 1.1
    \put{$g^\vee s^\vee_1$} at -0.5 0.2
    \put{$X^{\omega_1}$} at -0.5 1.45
    \put{$g^\vee s^\vee_0$} at .5 0.2
    \put{$X^{s_1\omega_1}$} at 1 0.65
    \put{$X^{w_0\omega_1}$} at -0.5 -0.2
    }%
    \put{\beginpicture
    \setplotsymbol({\tiny{$\bullet$}})
    \plot 1 0 0.5 0.866 /
    \plot 0.5 0.866 -0.5 0.866 /
    \plot -0.5 0.866 -1 0 /
    \plot -1 0 -0.5 -0.866 /
    \plot -0.5 -0.866 0.5 -0.866 /
    \plot 0.5 -0.866 1 0 /
    \put{$\bullet$} at 0 0
    \endpicture} at -0.5 0.866
    \put{\beginpicture
    \setplotsymbol({\tiny{$\bullet$}})
    \plot 1 0 0.5 0.866 /
    \plot 0.5 0.866 -0.5 0.866 /
    \plot -0.5 0.866 -1 0 /
    \plot -1 0 -0.5 -0.866 /
    \plot -0.5 -0.866 0.5 -0.866 /
    \plot 0.5 -0.866 1 0 /
    \put{$\bullet$} at 0 0
    \endpicture} at -0.5 -0.866
    \put{\beginpicture
    \setplotsymbol({\tiny{$\bullet$}})
    \plot 1 0 0.5 0.866 /
    \plot 0.5 0.866 -0.5 0.866 /
    \plot -0.5 0.866 -1 0 /
    \plot -1 0 -0.5 -0.866 /
    \plot -0.5 -0.866 0.5 -0.866 /
    \plot 0.5 -0.866 1 0 /
    \put{$\bullet$} at 0 0
    \endpicture} at -0.5 -2.598
    \put{\beginpicture
    \setplotsymbol({\tiny{$\bullet$}})
    \plot 1 0 0.5 0.866 /
    \plot 0.5 0.866 -0.5 0.866 /
    \plot -0.5 0.866 -1 0 /
    \plot -1 0 -0.5 -0.866 /
    \plot -0.5 -0.866 0.5 -0.866 /
    \plot 0.5 -0.866 1 0 /
    \put{$\bullet$} at 0 0
    \endpicture} at -0.5 2.598
    \put{\beginpicture
    \setplotsymbol({\tiny{$\bullet$}})
    \plot 1 0 0.5 0.866 /
    \plot 0.5 0.866 -0.5 0.866 /
    \plot -0.5 0.866 -1 0 /
    \plot -1 0 -0.5 -0.866 /
    \plot -0.5 -0.866 0.5 -0.866 /
    \plot 0.5 -0.866 1 0 /
    \put{$\bullet$} at 0 0
    \endpicture} at 1 0
    \put{\beginpicture
    \setplotsymbol({\tiny{$\bullet$}})
    \plot 1 0 0.5 0.866 /
    \plot 0.5 0.866 -0.5 0.866 /
    \plot -0.5 0.866 -1 0 /
    \plot -1 0 -0.5 -0.866 /
    \plot -0.5 -0.866 0.5 -0.866 /
    \plot 0.5 -0.866 1 0 /
    \put{$\bullet$} at 0 0
    \endpicture} at 1 1.732
    \put{\beginpicture
    \setplotsymbol({\tiny{$\bullet$}})
    \plot 1 0 0.5 0.866 /
    \plot 0.5 0.866 -0.5 0.866 /
    \plot -0.5 0.866 -1 0 /
    \plot -1 0 -0.5 -0.866 /
    \plot -0.5 -0.866 0.5 -0.866 /
    \plot 0.5 -0.866 1 0 /
    \put{$\bullet$} at 0 0
    \endpicture} at 1 -1.732
    \put{\beginpicture
    \setplotsymbol({\tiny{$\bullet$}})
    \plot 1 0 0.5 0.866 /
    \plot 0.5 0.866 -0.5 0.866 /
    \plot -0.5 0.866 -1 0 /
    \plot -1 0 -0.5 -0.866 /
    \plot -0.5 -0.866 0.5 -0.866 /
    \plot 0.5 -0.866 1 0 /
    \put{$\bullet$} at 0 0
    \endpicture} at 2.5 0.866
    \put{\beginpicture
    \setplotsymbol({\tiny{$\bullet$}})
    \plot 1 0 0.5 0.866 /
    \plot 0.5 0.866 -0.5 0.866 /
    \plot -0.5 0.866 -1 0 /
    \plot -1 0 -0.5 -0.866 /
    \plot -0.5 -0.866 0.5 -0.866 /
    \plot 0.5 -0.866 1 0 /
    \put{$\bullet$} at 0 0
    \endpicture} at 2.5 -0.866
    \put{\beginpicture
    \setplotsymbol({\tiny{$\bullet$}})
    \plot 1 0 0.5 0.866 /
    \plot 0.5 0.866 -0.5 0.866 /
    \plot -0.5 0.866 -1 0 /
    \plot -1 0 -0.5 -0.866 /
    \plot -0.5 -0.866 0.5 -0.866 /
    \plot 0.5 -0.866 1 0 /
    \put{$\bullet$} at 0 0
    \endpicture} at 2.5 -2.598
    \put{\beginpicture
    \setplotsymbol({\tiny{$\bullet$}})
    \plot 1 0 0.5 0.866 /
    \plot 0.5 0.866 -0.5 0.866 /
    \plot -0.5 0.866 -1 0 /
    \plot -1 0 -0.5 -0.866 /
    \plot -0.5 -0.866 0.5 -0.866 /
    \plot 0.5 -0.866 1 0 /
    \put{$\bullet$} at 0 0
    \endpicture} at 2.5 2.598
    \put{\beginpicture
    \setplotsymbol({\tiny{$\bullet$}})
    \plot 1 0 0.5 0.866 /
    \plot 0.5 0.866 -0.5 0.866 /
    \plot -0.5 0.866 -1 0 /
    \plot -1 0 -0.5 -0.866 /
    \plot -0.5 -0.866 0.5 -0.866 /
    \plot 0.5 -0.866 1 0 /
    \put{$\bullet$} at 0 0
    \endpicture} at -2 0
    \put{\beginpicture
    \setplotsymbol({\tiny{$\bullet$}})
    \plot 1 0 0.5 0.866 /
    \plot 0.5 0.866 -0.5 0.866 /
    \plot -0.5 0.866 -1 0 /
    \plot -1 0 -0.5 -0.866 /
    \plot -0.5 -0.866 0.5 -0.866 /
    \plot 0.5 -0.866 1 0 /
    \put{$\bullet$} at 0 0
    \endpicture} at -2 1.732
    \put{\beginpicture
    \setplotsymbol({\tiny{$\bullet$}})
    \plot 1 0 0.5 0.866 /
    \plot 0.5 0.866 -0.5 0.866 /
    \plot -0.5 0.866 -1 0 /
    \plot -1 0 -0.5 -0.866 /
    \plot -0.5 -0.866 0.5 -0.866 /
    \plot 0.5 -0.866 1 0 /
    \put{$\bullet$} at 0 0
    \endpicture} at -2 -1.732
    \plot 2.5 -4.33 3.7 -2.2516 /
    \plot 1.5 -4.33 3.7 -0.52 /
    \plot 0.5 -4.33 3.7 1.212 /
    \plot -0.5 -4.33 3.7 3.044 /
    \plot -1.5 -4.33 3.2 3.81 /
    \plot -2.5 -4.33 2.2 3.81 /
    \plot -3.5 -4.33 1.2 3.81 /
    \plot -3.7 -3.044 0.2 3.81 /
    \plot -3.7 -1.212 -0.8 3.81 /
    \plot -3.7 0.520 -1.8 3.81 /
    \plot -3.7 2.2516 -2.8 3.81 /
    \plot 2.5 4.33 3.7 2.2516 /
    \plot 1.5 4.33 3.7 0.52 /
    \plot 0.5 4.33 3.7 -1.212 /
    \plot -0.5 4.33 3.7 -3.044 /
    \plot -1.5 4.33 3.2 -3.81 /
    \plot -2.5 4.33 2.2 -3.81 /
    \plot -3.5 4.33 1.2 -3.81 /
    \plot -3.7 3.044 0.2 -3.81 /
    \plot -3.7 1.212 -0.8 -3.81 /
    \plot -3.7 -0.520 -1.8 -3.81 /
    \plot -3.7 -2.2516 -2.8 -3.81 /
    \plot -3.7 -3.464 4.5 -3.464 /
    \plot -3.7 -2.598 4.5 -2.598 /
    \plot -3.7 -1.732 4.5 -1.732 /
    \plot -3.7 -0.866 4.5 -0.866 /
    \plot -3.7 0 4.5 0 /
    \plot -3.7 3.464 4.5 3.464 /
    \plot -3.7 2.598 4.5 2.598 /
    \plot -3.7 1.732 4.5 1.732 /
    \plot -3.7 0.866 4.5 0.866 /
    {\small\put{$\fh^{-\alpha_1^\vee+d}$} at -3.5 -4.7
    \put{$\fh^{\alpha_1^\vee}$} at -2.5 -4.7
    \put{$\fh^{\alpha_1^\vee+2d}$} at -0.5 -4.7
    \put{$\fh^{\alpha_1^\vee+4d}$} at 1.5 -4.7 } %
    {\small\put{$\fh^{\alpha_2^\vee+d}$} at -3.5 4.7
    \put{$\fh^{\alpha_2^\vee}$} at -2.5 4.7
    \put{$\fh^{-\alpha_2^\vee+2d}$} at -0.5 4.7
    \put{$\fh^{-\alpha_2^\vee+4d}$} at 1.5 4.7 } %
    {\small\put{$\fh^{-\varphi^\vee+4d}$}[l] at 4.5 3.464
    \put{$\fh^{-\varphi^\vee+3d}$}[l] at 4.5 2.598
    \put{$\fh^{-\varphi^\vee+2d}$}[l] at 4.5 1.732
    \put{$\fh^{\alpha_0^\vee}$}[l] at 4.5 0.866
    \put{$\fh^{\varphi^\vee}$}[l] at 4.5 0
    \put{$\fh^{\varphi^\vee+d}$}[l] at 4.5 -0.866
    \put{$\fh^{\varphi^\vee+2d}$}[l] at 4.5 -1.732
    \put{$\fh^{\varphi^\vee+3d}$}[l] at 4.5 -2.598
    \put{$\fh^{\varphi^\vee+4d}$}[l] at 4.5 -3.464 }%
    \put{$+$} at 4.25 0.2
    \put{$+$} at 4.25 1.066
    \put{$+$} at 4.25 1.932
    \put{$+$} at 4.25 2.798
    \put{$+$} at 4.25 3.664
    \put{$-$} at 4.25 -1.066
    \put{$-$} at 4.25 -1.932
    \put{$-$} at 4.25 -2.798
    \put{$-$} at 4.25 -3.664
    \put{$+$} at 4.25 -3.264
    \put{$+$} at 4.25 -2.398
    \put{$+$} at 4.25 -1.532
    \put{$+$} at 4.25 -0.666
    \put{$-$} at 4.25 3.264
    \put{$-$} at 4.25 2.398
    \put{$-$} at 4.25 1.532
    \put{$-$} at 4.25 0.666
    \put{$-$} at 4.25 -0.2
    \put{$+$} at -0.9 -4.2
    \put{$-$} at -0.4 -4.2
    \put{$+$} at -1.9 -4.2
    \put{$-$} at -1.4 -4.2
    \put{$+$} at -2.9 -4.2
    \put{$-$} at -2.4 -4.2
    \put{$+$} at -3.9 -4.2
    \put{$-$} at -3.4 -4.2
    \put{$+$} at 0.1 -4.2
    \put{$-$} at 0.6 -4.2
    \put{$+$} at 1.1 -4.2
    \put{$-$} at 1.6 -4.2
    \put{$+$} at 2.1 -4.2
    \put{$-$} at 2.6 -4.2
    \put{$-$} at -0.9 4.2
    \put{$+$} at -0.4 4.2
    \put{$-$} at -1.9 4.2
    \put{$+$} at -1.4 4.2
    \put{$-$} at -2.9 4.2
    \put{$+$} at -2.4 4.2
    \put{$-$} at -3.9 4.2
    \put{$+$} at -3.4 4.2
    \put{$-$} at 0.1 4.2
    \put{$+$} at 0.6 4.2
    \put{$-$} at 1.1 4.2
    \put{$+$} at 1.6 4.2
    \put{$-$} at 2.1 4.2
    \put{$+$} at 2.6 4.2
    \put{Sheet $g^\vee$} at 0 -5.3
\endpicture$$

$$\beginpicture
\setcoordinatesystem units <1.4cm,1.4cm>         
\setplotarea x from -5 to 5, y from -5.5 to 5.2  
    \plot 2.8 -3.81 3.7 -2.2516 /
    \plot 1.8 -3.81 3.7 -0.52 /
    \plot 0.8 -3.81 3.7 1.212 /
    \plot -0.2 -3.81 3.7 3.044 /
    \plot -1.2 -3.81 3.2 3.81 /
    \plot -2.2 -3.81 2.2 3.81 /
    \plot -3.2 -3.81 1.2 3.81 /
    \plot -3.7 -3.044 0.2 3.81 /
    \plot -3.7 -1.212 -0.8 3.81 /
    \plot -3.7 0.520 -1.8 3.81 /
    \plot -3.7 2.2516 -2.8 3.81 /
    \plot -2.8 -3.81 -3.7 -2.2516 /
    \plot -1.8 -3.81 -3.7 -0.52 /
    \plot -0.8 -3.81 -3.7 1.212 /
    \plot 0.2 -3.81 -3.7 3.044 /
    \plot 1.2 -3.81 -3.2 3.81 /
    \plot 2.2 -3.81 -2.2 3.81 /
    \plot 3.2 -3.81 -1.2 3.81 /
    \plot 3.7 -3.044 -0.2 3.81 /
    \plot 3.7 -1.212 0.8 3.81 /
    \plot 3.7 0.520 1.8 3.81 /
    \plot 3.7 2.2516 2.8 3.81 /
    \plot -3.7 -3.464 3.7 -3.464 /
    \plot -3.7 -2.598 3.7 -2.598 /
    \plot -3.7 -1.732 3.7 -1.732 /
    \plot -3.7 -0.866 3.7 -0.866 /
    \plot -3.7 0 3.7 0 /
    \plot -3.7 3.464 3.7 3.464 /
    \plot -3.7 2.598 3.7 2.598 /
    \plot -3.7 1.732 3.7 1.732 /
    \plot -3.7 0.866 3.7 0.866 /
    \put{$\bullet$} at 0 0
    \put{$\bullet$} at 2 3.464
    \put{$\bullet$} at 2 -3.464
    \put{$\bullet$} at -1 3.464
    \put{$\bullet$} at -1 -3.464
    \small{\put{$\kappa$} at 0 0.6
    \put{$\kappa s^\vee_0$} at -0.5 0.25
    \put{$\kappa s^\vee_1$} at 0.5 0.25
    \put{$\kappa s^\vee_2$} at 0 1.1
    \put{$X^{\omega_2}$} at 0.5 1.5
    \put{$X^{s_2\omega_2}$} at -1 0.7
    \put{$X^{w_0\omega_2}$} at 0.5 -0.17
    }%
    \put{\beginpicture
    \setplotsymbol({\tiny{$\bullet$}})
    \plot 1 0 0.5 0.866 /
    \plot 0.5 0.866 -0.5 0.866 /
    \plot -0.5 0.866 -1 0 /
    \plot -1 0 -0.5 -0.866 /
    \plot -0.5 -0.866 0.5 -0.866 /
    \plot 0.5 -0.866 1 0 /
    \put{$\bullet$} at 0 0
    \endpicture} at 0.5 0.866
    \put{\beginpicture
    \setplotsymbol({\tiny{$\bullet$}})
    \plot 1 0 0.5 0.866 /
    \plot 0.5 0.866 -0.5 0.866 /
    \plot -0.5 0.866 -1 0 /
    \plot -1 0 -0.5 -0.866 /
    \plot -0.5 -0.866 0.5 -0.866 /
    \plot 0.5 -0.866 1 0 /
    \put{$\bullet$} at 0 0
    \endpicture} at 0.5 -0.866
    \put{\beginpicture
    \setplotsymbol({\tiny{$\bullet$}})
    \plot 1 0 0.5 0.866 /
    \plot 0.5 0.866 -0.5 0.866 /
    \plot -0.5 0.866 -1 0 /
    \plot -1 0 -0.5 -0.866 /
    \plot -0.5 -0.866 0.5 -0.866 /
    \plot 0.5 -0.866 1 0 /
    \put{$\bullet$} at 0 0
    \endpicture} at 0.5 -2.598
    \put{\beginpicture
    \setplotsymbol({\tiny{$\bullet$}})
    \plot 1 0 0.5 0.866 /
    \plot 0.5 0.866 -0.5 0.866 /
    \plot -0.5 0.866 -1 0 /
    \plot -1 0 -0.5 -0.866 /
    \plot -0.5 -0.866 0.5 -0.866 /
    \plot 0.5 -0.866 1 0 /
    \put{$\bullet$} at 0 0
    \endpicture} at 0.5 2.598
    \put{\beginpicture
    \setplotsymbol({\tiny{$\bullet$}})
    \plot 1 0 0.5 0.866 /
    \plot 0.5 0.866 -0.5 0.866 /
    \plot -0.5 0.866 -1 0 /
    \plot -1 0 -0.5 -0.866 /
    \plot -0.5 -0.866 0.5 -0.866 /
    \plot 0.5 -0.866 1 0 /
    \put{$\bullet$} at 0 0
    \endpicture} at -1 0
    \put{\beginpicture
    \setplotsymbol({\tiny{$\bullet$}})
    \plot 1 0 0.5 0.866 /
    \plot 0.5 0.866 -0.5 0.866 /
    \plot -0.5 0.866 -1 0 /
    \plot -1 0 -0.5 -0.866 /
    \plot -0.5 -0.866 0.5 -0.866 /
    \plot 0.5 -0.866 1 0 /
    \put{$\bullet$} at 0 0
    \endpicture} at -1 1.732
    \put{\beginpicture
    \setplotsymbol({\tiny{$\bullet$}})
    \plot 1 0 0.5 0.866 /
    \plot 0.5 0.866 -0.5 0.866 /
    \plot -0.5 0.866 -1 0 /
    \plot -1 0 -0.5 -0.866 /
    \plot -0.5 -0.866 0.5 -0.866 /
    \plot 0.5 -0.866 1 0 /
    \put{$\bullet$} at 0 0
    \endpicture} at -1 -1.732
    \put{\beginpicture
    \setplotsymbol({\tiny{$\bullet$}})
    \plot 1 0 0.5 0.866 /
    \plot 0.5 0.866 -0.5 0.866 /
    \plot -0.5 0.866 -1 0 /
    \plot -1 0 -0.5 -0.866 /
    \plot -0.5 -0.866 0.5 -0.866 /
    \plot 0.5 -0.866 1 0 /
    \put{$\bullet$} at 0 0
    \endpicture} at -2.5 0.866
    \put{\beginpicture
    \setplotsymbol({\tiny{$\bullet$}})
    \plot 1 0 0.5 0.866 /
    \plot 0.5 0.866 -0.5 0.866 /
    \plot -0.5 0.866 -1 0 /
    \plot -1 0 -0.5 -0.866 /
    \plot -0.5 -0.866 0.5 -0.866 /
    \plot 0.5 -0.866 1 0 /
    \put{$\bullet$} at 0 0
    \endpicture} at -2.5 -0.866
    \put{\beginpicture
    \setplotsymbol({\tiny{$\bullet$}})
    \plot 1 0 0.5 0.866 /
    \plot 0.5 0.866 -0.5 0.866 /
    \plot -0.5 0.866 -1 0 /
    \plot -1 0 -0.5 -0.866 /
    \plot -0.5 -0.866 0.5 -0.866 /
    \plot 0.5 -0.866 1 0 /
    \put{$\bullet$} at 0 0
    \endpicture} at -2.5 -2.598
    \put{\beginpicture
    \setplotsymbol({\tiny{$\bullet$}})
    \plot 1 0 0.5 0.866 /
    \plot 0.5 0.866 -0.5 0.866 /
    \plot -0.5 0.866 -1 0 /
    \plot -1 0 -0.5 -0.866 /
    \plot -0.5 -0.866 0.5 -0.866 /
    \plot 0.5 -0.866 1 0 /
    \put{$\bullet$} at 0 0
    \endpicture} at -2.5 2.598
    \put{\beginpicture
    \setplotsymbol({\tiny{$\bullet$}})
    \plot 1 0 0.5 0.866 /
    \plot 0.5 0.866 -0.5 0.866 /
    \plot -0.5 0.866 -1 0 /
    \plot -1 0 -0.5 -0.866 /
    \plot -0.5 -0.866 0.5 -0.866 /
    \plot 0.5 -0.866 1 0 /
    \put{$\bullet$} at 0 0
    \endpicture} at 2 0
    \put{\beginpicture
    \setplotsymbol({\tiny{$\bullet$}})
    \plot 1 0 0.5 0.866 /
    \plot 0.5 0.866 -0.5 0.866 /
    \plot -0.5 0.866 -1 0 /
    \plot -1 0 -0.5 -0.866 /
    \plot -0.5 -0.866 0.5 -0.866 /
    \plot 0.5 -0.866 1 0 /
    \put{$\bullet$} at 0 0
    \endpicture} at 2 1.732
    \put{\beginpicture
    \setplotsymbol({\tiny{$\bullet$}})
    \plot 1 0 0.5 0.866 /
    \plot 0.5 0.866 -0.5 0.866 /
    \plot -0.5 0.866 -1 0 /
    \plot -1 0 -0.5 -0.866 /
    \plot -0.5 -0.866 0.5 -0.866 /
    \plot 0.5 -0.866 1 0 /
    \put{$\bullet$} at 0 0
    \endpicture} at 2 -1.732
    \plot 2.5 -4.33 3.7 -2.2516 /
    \plot 1.5 -4.33 3.7 -0.52 /
    \plot 0.5 -4.33 3.7 1.212 /
    \plot -0.5 -4.33 3.7 3.044 /
    \plot -1.5 -4.33 3.2 3.81 /
    \plot -2.5 -4.33 2.2 3.81 /
    \plot -3.5 -4.33 1.2 3.81 /
    \plot -3.7 -3.044 0.2 3.81 /
    \plot -3.7 -1.212 -0.8 3.81 /
    \plot -3.7 0.520 -1.8 3.81 /
    \plot -3.7 2.2516 -2.8 3.81 /
    \plot 2.5 4.33 3.7 2.2516 /
    \plot 1.5 4.33 3.7 0.52 /
    \plot 0.5 4.33 3.7 -1.212 /
    \plot -0.5 4.33 3.7 -3.044 /
    \plot -1.5 4.33 3.2 -3.81 /
    \plot -2.5 4.33 2.2 -3.81 /
    \plot -3.5 4.33 1.2 -3.81 /
    \plot -3.7 3.044 0.2 -3.81 /
    \plot -3.7 1.212 -0.8 -3.81 /
    \plot -3.7 -0.520 -1.8 -3.81 /
    \plot -3.7 -2.2516 -2.8 -3.81 /
    \plot -3.7 -3.464 4.5 -3.464 /
    \plot -3.7 -2.598 4.5 -2.598 /
    \plot -3.7 -1.732 4.5 -1.732 /
    \plot -3.7 -0.866 4.5 -0.866 /
    \plot -3.7 0 4.5 0 /
    \plot -3.7 3.464 4.5 3.464 /
    \plot -3.7 2.598 4.5 2.598 /
    \plot -3.7 1.732 4.5 1.732 /
    \plot -3.7 0.866 4.5 0.866 /
    {\small\put{$\fh^{-\alpha_1^\vee+d}$} at -3.5 -4.7
    \put{$\fh^{\alpha_1^\vee}$} at -2.5 -4.7
    \put{$\fh^{\alpha_1^\vee+2d}$} at -0.5 -4.7
    \put{$\fh^{\alpha_1^\vee+4d}$} at 1.5 -4.7 } %
    {\small\put{$\fh^{\alpha_2^\vee+d}$} at -3.5 4.7
    \put{$\fh^{\alpha_2^\vee}$} at -2.5 4.7
    \put{$\fh^{-\alpha_2^\vee+2d}$} at -0.5 4.7
    \put{$\fh^{-\alpha_2^\vee+4d}$} at 1.5 4.7 } %
    {\small\put{$\fh^{-\varphi^\vee+4d}$}[l] at 4.5 3.464
    \put{$\fh^{-\varphi^\vee+3d}$}[l] at 4.5 2.598
    \put{$\fh^{-\varphi^\vee+2d}$}[l] at 4.5 1.732
    \put{$\fh^{\alpha_0^\vee}$}[l] at 4.5 0.866
    \put{$\fh^{\varphi^\vee}$}[l] at 4.5 0
    \put{$\fh^{\varphi^\vee+d}$}[l] at 4.5 -0.866
    \put{$\fh^{\varphi^\vee+2d}$}[l] at 4.5 -1.732
    \put{$\fh^{\varphi^\vee+3d}$}[l] at 4.5 -2.598
    \put{$\fh^{\varphi^\vee+4d}$}[l] at 4.5 -3.464 }%
    \put{$+$} at 4.25 0.2
    \put{$+$} at 4.25 1.066
    \put{$+$} at 4.25 1.932
    \put{$+$} at 4.25 2.798
    \put{$+$} at 4.25 3.664
    \put{$-$} at 4.25 -1.066
    \put{$-$} at 4.25 -1.932
    \put{$-$} at 4.25 -2.798
    \put{$-$} at 4.25 -3.664
    \put{$+$} at 4.25 -3.264
    \put{$+$} at 4.25 -2.398
    \put{$+$} at 4.25 -1.532
    \put{$+$} at 4.25 -0.666
    \put{$-$} at 4.25 3.264
    \put{$-$} at 4.25 2.398
    \put{$-$} at 4.25 1.532
    \put{$-$} at 4.25 0.666
    \put{$-$} at 4.25 -0.2
    \put{$+$} at -0.9 -4.2
    \put{$-$} at -0.4 -4.2
    \put{$+$} at -1.9 -4.2
    \put{$-$} at -1.4 -4.2
    \put{$+$} at -2.9 -4.2
    \put{$-$} at -2.4 -4.2
    \put{$+$} at -3.9 -4.2
    \put{$-$} at -3.4 -4.2
    \put{$+$} at 0.1 -4.2
    \put{$-$} at 0.6 -4.2
    \put{$+$} at 1.1 -4.2
    \put{$-$} at 1.6 -4.2
    \put{$+$} at 2.1 -4.2
    \put{$-$} at 2.6 -4.2
    \put{$-$} at -0.9 4.2
    \put{$+$} at -0.4 4.2
    \put{$-$} at -1.9 4.2
    \put{$+$} at -1.4 4.2
    \put{$-$} at -2.9 4.2
    \put{$+$} at -2.4 4.2
    \put{$-$} at -3.9 4.2
    \put{$+$} at -3.4 4.2
    \put{$-$} at 0.1 4.2
    \put{$+$} at 0.6 4.2
    \put{$-$} at 1.1 4.2
    \put{$+$} at 1.6 4.2
    \put{$-$} at 2.1 4.2
    \put{$+$} at 2.6 4.2
    \put{Sheet $\kappa = (g^\vee)^2$} at 0 -5.3
\endpicture$$

\end{section}

\end{document}